%
%

%
%


\magnification 1200
\input amstex
\documentstyle{amsppt}
\NoBlackBoxes
\NoRunningHeads

\vsize = 9.4 truein

\define\ggm{{G}/\Gamma}

\define\vre{\varepsilon}
\define\hd{Hausdorff dimension}
\define\hs{homogeneous space}
\define\df{\overset\text{def}\to=}
\define\un#1#2{\underset\text{#1}\to#2}
\define\br{\Bbb R}
\define\bn{\Bbb N}
\define\bz{\Bbb Z}

\define\di{Diophantine}
\define\da{Diophantine approximation}
\define\de{Diophantine exponent}
\define\va{\bold a}
\define\vb{\bold b}
\define\ve{\bold e}
\define\vx{\bold x}
\define\vy{\bold y}
\define\vz{\bold z}
\define\vv{\bold v}

\define\vw{\bold w}
\define\vp{\bold p}
\define\vq{\bold q}

\define\vc{\bold c}
\define\vf{\bold f}
\define\vg{\bold g}
\define\vh{\bold h}
\define\vt{\bold t}

\define\nz{\smallsetminus \{0\}}

\define\cag{$(C,\alpha)$-good}

\define\vrn{\varnothing}
\define\ssm{\smallsetminus}
\define\vwa{very well approximable}

\define\mr{M_{m,n}}
\define\amr{$A\in M_{m,n}(\br)$}

\define\lio{\Lambda \in \Omega}
\def\ma{multiplicative approximation}
\define\rk{\operatorname{rk}}
\define\GL{\operatorname{GL}}
\define\SL{\operatorname{SL}}

\topmatter
\title 
Extremal subspaces and their submanifolds
\endtitle  

\author { Dmitry Kleinbock} \\ 
  { \rm 
   Brandeis University} 
\endauthor

    \address{ Dmitry Kleinbock,  Department of
Mathematics, Brandeis University, Waltham, MA 02454-9110}
  \endaddress

\email kleinboc\@brandeis.edu \endemail

  \thanks  
Supported in part by NSF
Grant DMS-0072565.
\endthanks

\abstract 
It was proved in the  paper \cite{KM1} that the properties of almost
all points of $\br^n$ being not very well (multiplicatively) approximable are inherited by 
nondegenerate in $\br^n$ (read: not contained in a proper affine subspace) smooth
submanifolds. In this paper we consider submanifolds which are
contained in proper affine subspaces, and prove that  the aforementioned
\di\ properties pass from a subspace to its nondegenerate
submanifold.
The  proofs are based on 
a correspondence between multidimensional \da\  and dynamics of lattices in 
Euclidean spaces.

 \endabstract

\dedicatory {To appear in Geom.~Funct.~Anal.} \enddedicatory


\endtopmatter
\document

\heading{1. Introduction}
\endheading 

We  denote by $\mr$ the space of real matrices with  $m$ rows and $n$
columns. $I_k \in M_{k,k}$ stands for the identity matrix.
Vectors are named by  lowercase boldface letters, such as $\vx =
(x_i\mid 1\le i \le k)$.
For $\vx\in\br^k$ we let 
$\|\vx\| =
\max_{1\le i
\le k}|x_i|
$. $0$ stands for a zero vector in any dimension, as
well as  a zero matrix of any size.  The Lebesgue measure in $\br^k$ will be denoted by
$|\cdot|$.

We start by recalling several basic facts from the theory of 
 \da.  
For $v > 0$ and $m,n\in\bn$, let us denote by $\Cal W_v(m,n)$ the set of 
matrices $A\in\mr$ for which  there are infinitely
many
$\vq\in
\bz^n$ such that
$$
 \|A\vq + \vp\|   \le \|\vq\|^{-v} \quad \text{for some
}\vp\in\bz^m\,.  \tag 1.1
$$
Clearly $\Cal
W_{v_1}(m,n) \supset \Cal
W_{v_2}(m,n)$ if $v_1 \le v_2$. We will also use the notation
$$
\Cal W_v^{\sssize +}(m,n)\df \bigcup_{u > v} \Cal W_{u}(m,n)\quad\text{and}\quad\Cal W_v^{\sssize -}(m,n)\df
\bigcap_{u < v} \Cal W_{u}(m,n)\,.
$$

One knows that $\Cal W_{n/m}(m,n) = \mr$ by Dirichlet's
Theorem, and that the Lebesgue measure of $\Cal W_v(m,n)$ is zero whenever $ v > n/m$ due to the
Borel-Cantelli Lemma. In particular, the set  $\Cal W_{n/m}^{\sssize +}(m,n)$ has zero measure. Matrices
from the latter set are called {\sl \vwa\/}, to be abbreviated as VWA.

It follows from Khintchine's
Transference Principle, see e.g.~\cite{C, Chapter V}, that the
statements $A\in \Cal W_{n/m}^{\sssize +}(m,n)$  and $A^{\sssize T}\in \Cal W_{m/n}^{\sssize +}(n,m)$ are
equivalent. In particular, a vector $\vy\in\br^n$ interpreted as an
$n\times 1$ matrix is VWA iff it is VWA when viewed as a $1\times n$
matrix. Our goal in the present paper is to look at VWA vectors in $\br^n$, and 
it will be more convenient for us to use the {\it row vector\/}
approach, so that, for $\vy$ as above and for $\vq\in\br^n$, $\vy\vq$ will stand for 
$y_1q_1 + \dots + y_nq_n$.  In
view of the aforementioned duality, this causes no loss of generality, and, hopefully, will
cause no confusion.

\medskip

We now specialize to the case $m = 1$; that is, consider \di\ properties
of vectors ($=$ row matrices) $\vy\in\br^n$. Following the 
terminology introduced by V.~Sprind\v zuk, say that a 
submanifold 
 $\Cal M$ of $\br^n$ is  {\sl extremal\/} if almost all $\vy\in \Cal M$ (with
respect to the natural measure class) are not VWA. In other words, 
if the property of  generic $\vy\in\br^n$ being not VWA is
{\it inherited\/} by the submanifold. Pushing this terminology a little further, 
let us say that a 
map $\vf$ from an open subset $U$ of $\br^d$ to $\br^n$  is 
{\sl extremal\/} if $\vf(\vx)$ is not VWA for a.e.~$\vx\in
U$. 

Proving
extremality of smooth manifolds/maps has been one of the central topics of metric \da\ for the
last 40 years, the major driving force being 
Sprind\v zuk's 1964 solution \cite{Sp1} of a long-standing problem of K.~Mahler \cite{M}, that
is,  proving the extremality of the so-called {\sl rational
normal\/} or {\sl Veronese} curve  
 $$
\Cal M = \{(x,x^2,\dots,x^n)\mid x\in \br\}\,.\tag 1.2
$$ 
See \cite{Sp2, Sp3, BD} for history and
references. 


In his 1980 survey of the field \cite{Sp4}, Sprind\v zuk conjectured that a
real analytic manifold $\Cal M$ is extremal whenever it  is not contained in any proper affine
subspace of $\br^n$. The latter condition, loosely put, says
that $\Cal M$ `remembers'  the dimension of the space it is imbedded into; and
the conjecture asserts that $\Cal M$ must also `remember' the law of almost
all points being not VWA.

This conjecture was proved by G.\,A.~Margulis and the author \cite{KM1} in a
stronger form, with the aforementioned geometric condition replaced by an
analytic one, and the real analytic class extended to $C^k$ for large enough
$k$. We need the following definitions. Let $U$ be an open subset of
 $\br^d$,  $\Cal L$  an affine
subspace of $\br^n$, and let  $\vf = (f_1,\dots,f_n)$ be a
  $C^k$ 
map $U\to \Cal L$. For $l\le k$ and $\vx\in U$, say that 
$\vf$  is {\sl $l$-nondegenerate in $\Cal L$ at $\vx$\/} if 
$$
\aligned
\text{the linear part of  $\Cal L$ is spanned by }\\ \text{partial derivatives of $\vf$ at
$\vx$ of order up}&\text{ to $l$}
\endaligned
\tag 1.3
$$
 (a linear subspace $\Cal L_0$ of  $\br^n$ is called the {\sl linear part\/} of $\Cal L$ if 
$\Cal L = \Cal L_0 + \vy$
for some $\vy \in\br^n$). We will say that
$\vf$ is 
 {\sl nondegenerate in $\Cal L$ at $\vx$\/}  if (1.3) holds for
some $l$. 
If $\Cal M$ is a $d$-dimensional 
submanifold of $\Cal L$, we will say
that $\Cal M$ is {\sl nondegenerate  in $\Cal L$ at $\vy\in \Cal M$} if any
(equivalently, some)  diffeomorphism
$\vf$ between an open subset $U$  of $\br^d$ and a neighborhood of
$\vy$ in $\Cal M$ is   nondegenerate  in $\Cal L$ at $\vf^{-1}(\vy)$. We will say
that
$\vf:U\to \Cal L$ (resp.~$\Cal M\subset \Cal L$) is {\sl nondegenerate  in
$\Cal L$\/}  if it is nondegenerate  in $\Cal L$ at almost every point of $U$ 
(resp.~$\Cal M$, in the sense of the 
natural measure class on $\Cal M$). 

One of the main results of \cite{KM1} is the following 

\proclaim{Theorem 1.1} Let $\vf:U\to\br^n$, $U\subset\br^d$, be a  smooth 
map which is nondegenerate  in $\br^n$. Then $\vf$ is
extremal.
\endproclaim

In particular, smooth submanifolds of $\br^n$ which are nondegenerate  in $\br^n$ are
extremal. Note that many special cases were proved before  the general case; see
\cite{KM1, BD} for a detailed account, and \cite{BKM, BBKM, KLW} for further developments.

The goal of the present paper is to study manifolds for which the
aforementioned non-degeneracy-in-$\br^n$ condition fails. In fact, the
simplest ones, namely proper affine subspaces of $\br^n$ themselves, have
been the subject of several papers \cite{S2, Sp3, BBDD}, and certain
conditions have been found sufficient for their extremality. To the best
of the author's knowledge, nobody has yet turned attention to proper
submanifolds of affine subspaces of $\br^n$. Let us now state one  of the main
results of the present paper, which addresses this gap.

\proclaim{Theorem 1.2} Let $\Cal L$ be an 
affine subspace  of $\br^n$. Then:
\roster
\item"{\rm (a)}" if $\Cal L$  is extremal and $\vf:U\to \Cal L$, $U\subset\br^d$,  is  a  smooth 
map  which is nondegenerate  in $\Cal L$, then $\vf$ is extremal;
\item"{\rm (b)}" if $\Cal L$  is not extremal, then all points of $\Cal L$ 
are VWA (in particular, no subset of $\Cal L$  is extremal).
\endroster
\endproclaim

This result generalizes Theorem 1.1, showing that the extremality of affine subspaces
is inherited by their nondegenerate submanifolds.  It also implies that a manifold
nondegenerate in some affine subspace 
of
$\br^n$ is extremal if and only if this subspace contains at least one not \vwa\ point. 
Cf.~a similar statement conjectured by B.~Weiss in the context of  interval exchange
transformations and Teichm\"uller flows
\cite{W, Conjecture 2.1}.

The  proof is based on the methods of  \cite{KM1}, that is, on the
correspondence between approximation properties of vectors and trajectories
of  lattices in Euclidean spaces. Necessary background is reviewed in \S 2. 
Then in  \S 3  we use the language
of lattices to give a necessary and sufficient   condition for a map $\vf:U\to\br^n$,
$U\subset\br^d$,  within a certain class of {\it good\/} maps to be extremal, and then
show that  this condition is inherited by nondegenerate submanifolds of affine subspaces. 

Dealing with an ${s}$-dimensional affine subspace  of $\br^n$, one can
be more specific and phrase the aforementioned condition in terms of coefficients of an affine
map parametrizing the subspace. By permuting variables one can  without loss of generality 
choose 
a parametrizing map  of the form
$\vx\mapsto(\vx,
\vx A' + \va_0)$,
where 
$A'$ is a matrix of size ${s} \times(n-{s})$ and $\va_0\in\br^{n-{s}}$ (here both $\vx$ and $\va_0$ are
row vectors). In an even more abbreviated way, we will denote   the vector $(1,x_1,\dots,x_{s})$ by
$\tilde
\vx$, and the matrix
$\left(\matrix
\va_0 \\ A'
\endmatrix\right)$ by $A\in M_{{s}+1,{n-{s}}}$; then $\Cal L$ is parametrized by $$
\vx\mapsto(\vx,
\tilde\vx A)\,.\tag 1.4
$$
We show in \S 4 how the results of \S 3 allow one to write down a  condition on $A$
(see Theorem 4.3)
equivalent to the extremality of the map (1.4).  On the other hand, 
it easily follows from the definitions, as explained in  \S 4, that  every point of
$\Cal L$ parametrized by (1.4) is VWA whenever 
$A$ belongs to
$\Cal W_{n}^{\sssize +}({s}+1,{n-{s}})$. 
We  show 
that
the converse is also true
 in the following two cases, and by the following two methods:
\pagebreak
$$
{s} = n-1 \text{ (that is, $\Cal L$ is an affine hyperplane)}\,,\tag 1.5
$$ 
--- as a consequence of
Theorem 4.3, 
and  
$$
\aligned
{s} = 1\text{  and $A$ is of the form $\left(\matrix
0 \\ \vb
\endmatrix\right)$ for a row vector $\vb$}\\ \text{ (that is, $\Cal L$ is a line passing through the
origin)}\,.\quad
\endaligned
\tag 1.6
$$ 
--- using an argument borrowed from  \cite{BBDD}. In other words, the following can be proved:

\proclaim{Theorem 1.3} In the two special cases {\rm (1.5)} and {\rm (1.6)}, the map
{\rm (1.4)} is extremal if and only if 
$$A\notin
\Cal W_{n}^{\sssize +}({s}+1,{n-{s}})\,.\tag 1.7$$ \endproclaim

Whether the same is true for an arbitrary affine subspace is
not clear. Since matrices $A$ as above provide local coordinate charts to the set  of 
${s}$-dimensional affine subspaces  of $\br^n$, and in view of
M.~Dodson's
\cite{Do} formula for the
\hd\ of the sets
$\Cal W_{v}(m,n)$, the affirmative answer to the above  question would imply that the dimension
of the (null) set of non-extremal
${s}$-dimensional affine subspaces  of $\br^n$ is equal to 
$$
\dim\big(\Cal W_{n}^{\sssize +}({s}+1,{n-{s}})\big)  = (n-{s}-1)({s}+1) + 1\,,\tag 1.8
$$
which is  precisely
$1\  +\ $ the \hd\ of the set of `rational' ${s}$-dimensional subspaces, i.e.~of the
set 
$$\big\{A\in M_{{s}+1,n-{s}}\bigm| A\vq +
\vp = 0\,\text{ for some }\vp\in\bz^{s+1},\vq\in\bz^{n-s}\nz\big\}\,.$$

Other
 open problems and generalizations are  discussed in the last two sections of the
paper. This includes the so-called {\it multiplicative\/} modification of the
standard set-up, which is the subject of \S 5. Namely, there  we define {\it not very well
multiplicatively approximable\/} (not VWMA, a property stronger than `not VWA' but still
generic in $\br^n$) vectors 
  and {\it strongly extremal\/} manifolds (i.e.~those for which almost all points
are not VWMA). It was proved in \cite{KM1} that smooth nondegenerate submanifolds of $\br^n$
 are strongly
extremal; we generalize this as follows:

\proclaim{Theorem 1.4} Let $\Cal L$ be an 
affine subspace  of $\br^n$. Then:
\roster
\item"{\rm (a)}" if $\Cal L$  is  strongly extremal and $\vf:U\to \Cal L$ is  a  smooth 
map  which is nondegenerate  in $\Cal L$, then $\vf$ is strongly  extremal;
\item"{\rm (b)}" if $\Cal L$  is not  strongly extremal, then all points of $\Cal L$ 
are VWMA (in particular, no subset of $\Cal L$  is  strongly extremal).
\endroster
\endproclaim

Similarly to Theorem 1.2, this is done by writing down a  necessary and sufficient  
condition (see Theorem 5.3) for a good map  to be strongly extremal, and then showing that  
condition to be inherited by nondegenerate submanifolds of affine subspaces. Following the lines
of \S 4, we are able to simplify that condition in the case (1.5), thus explicitly
describing strongly extremal hyperplanes and identifying those which are extremal but not strongly
extremal. Whether  this can be extended beyond the codimension one case is an open question.

\heading{2. \da\ and lattices}
\endheading

In this section we introduce some notation and terminology 
which will help us  work with
discrete subgroups $\Gamma$ of $\br^k$, $k\in \bn$. We 
define the {\sl rank\/} $\rk(\Gamma)$ 
of
$\Gamma$ to be the
dimension of $\br\Gamma$. 
Also define 
$\delta(\Gamma)$ to be the norm of a nonzero element of $\Gamma$ with the
smallest norm, that is, 
$$
\delta(\Gamma) \df \inf_{\vv\in\Gamma\nz}\|\vv\|\,.
$$
For $0\le j\le k$, let us denote by  $\Cal S_{k,j}$ the set of
all subgroups of $\bz^k$ of rank $j$, and by $\Cal
S_k$ the set of all
 nonzero   subgroups of $\bz^k$ of rank smaller than $k$, that is, 
$\Cal
S_k \df \cup_{j = 1}^{k-1} \Cal
S_{k,j}$. 

It  will be useful to consider 
exterior products
of vectors generating $\Gamma$.
Namely, if $\Gamma\in \Cal S_{k,j}$, say 
that $\vw\in \bigwedge^j(\br^{k})$ {\sl represents\/} $\Gamma$ if 
$$
\vw =  \cases &1\qquad\qquad\quad
\text{\ \ \ if } j = 0\\ &\vv_{1}\wedge\dots\wedge \vv_{j}\quad
\text{if }j > 0 \text{ and }\vv_{1},\dots, \vv_{j} \text{ is a basis of }
\Gamma\,.\endcases
$$
Clearly the element representing $\Gamma$ is defined up to a sign. 
With some abuse of notation, we will also denote by $\Cal
S_{k,j}$ and $\Cal
S_k$ the set of $\vw\in \bigwedge(\br^{k})$ representing $\Gamma\in\Cal
S_{k,j}$ and $\in\Cal
S_k$ respectively. 

The set of all 
{\sl lattices\/} (discrete subgroups   of maximal rank) in $\br^{k}$ of covolume one can be
identified with the \hs\
$\SL_{k}(\br)/\SL_{k}(\bz)$, which we will denote 
by $\Omega_k$. It is a noncompact space with
finite  $\SL_{k}(\br)$-invariant measure, and the restriction of
the function
$\delta(\cdot)$ defined above to this space can be used to describe
its geometry at infinity. Namely, Mahler's Compactness Criterion 
\cite{R, Corollary 10.9} says that a subset of $\Omega_k$ is relatively
compact if and only if $\delta$ is bounded away from zero on this subset.
Further, it follows from the {\it reduction
theory\/} for $\SL_{k}(\bz)$, see e.g.~\cite{Si, Satz 4}, that 
the ratio of $1 + \log\big(1/{\delta(\cdot)}\big)$ 
and $1 + \text{dist}(\cdot, \bz^{k})$ is bounded between two 
positive constants for any right invariant Riemannian metric 
`dist' on the space of lattices. In other words, a lattice $\lio_k$ for
which $\delta(\Lambda)$ is small is approximately 
$\log\big(1/{\delta(\Lambda)}\big)$ away from the base point 
$\bz^{k}$. The reader is referred to \cite{K1} for more details.

This justifies the following definition: for $\gamma \ge 0$ and
any one-parameter semigroup $
F = \{g_t\mid t\ge 0\}$ acting on $\Omega_k$, say that the $F$-trajectory of
$\Lambda\in
\Omega_k$  {\sl grows with exponent\/} 
$\ge \gamma$ if there exist arbitrarily large positive
$t$ such that 
$$
\delta( g_{t} \Lambda) \le e^{-\gamma t}\,.
$$
 Also define the {\sl growth exponent\/}
$\gamma_{\sssize F}(\Lambda)$ of $\Lambda$ with respect to $F$
to be the supremum of all $\gamma$ for which the $F$-trajectory of
$\Lambda$  grows with exponent 
$\ge \gamma$. In view of the preceding remark, one has
$$
\gamma_{\sssize F}(\Lambda) = \limsup_{t\to\infty}
\frac{\text{dist}(g_t\Lambda, \bz^{m})}t\,.
$$

Now let us describe a  correspondence, dating back to \cite{S3} and \cite{D}, 
between approximation properties of vectors  $\vy\in\br^n$ and dynamics of
certain trajectories  in $\Omega_{n+1}$. Given a
row vector  $\vy\in\br^n$ one  considers
 a lattice  $u_\vy\bz^{n+1}$ in $\br^{n+1}$, where 
$
u_\vy \df \left(\matrix
1 & \vy  \\
0 & I_n
\endmatrix \right)
$; that is, the collection of vectors  of the form 
$\left(\matrix \vy\vq + p \\  \vq \endmatrix
\right)
$,
where $p\in\bz$ and 
$\vq \in \bz^n$. Then one reads \di\ properties of
$\vy$  from the behavior of the trajectory 
$Fu_\vy\bz^{n+1}$,
where
$$
F = \{g_t\mid t\ge 0\}\,,\quad\text{with}\quad g_t = \text{\rm diag}(e^{t},
e^{-t/n},\dots,e^{-t/n})\,,\tag 2.1
$$
is  a one-parameter subsemigroup of
$\SL_{n+1}(\br)$  which  expands the first  coordinate and uniformly contracts the last $n$
coordinates of vectors in $\br^{n+1}$.  

\comment
Let us 
briefly discuss the geometry  of the space 
where the action will be taking place. 
 The set of all 
lattices in $\br^{n+1}$ of covolume one can be identified with the \hs\
$\SL_{n+1}(\br)/\SL_{n+1}(\bz)$, which we will denote 
by $\Omega$. It is a noncompact space with
finite  $\SL_{n+1}(\br)$-invariant measure, and the restriction of
the function
$\delta(\cdot)$ defined above to this space can be used to describe
its geometry at infinity. Namely, Mahler's Compactness Criterion 
\cite{R, Corollary 10.9} says that a subset of $\Omega$ is relatively
compact if and only if $\delta$ is bounded away from zero on this subset.
Further, it follows from the {\it reduction
theory\/} for arithmetic groups, see e.g.~\cite{Si, section 10}, that 
the ratio of $1 + \log\big(1/{\delta(\cdot)}\big)$ 
and $1 + \text{dist}(\cdot, \bz^{n+1})$ is bounded between two 
positive constants for any right invariant Riemannian metric 
`dist' on the space of lattices. In other words, a lattice $\lio$ for
which $\delta(\Lambda)$ is small is approximately 
$\log\big(1/{\delta(\Lambda)}\big)$ away from the base point 
$\bz^{n+1}\in\Omega$. The reader is referred to \cite{K} for more details.

This justifies the following definition: for $\gamma \ge 0$ and
any one-parameter semigroup $
g_t$, $t\ge 0$, acting on $\Omega$, say that the $g_{t}$-trajectory of
$\Lambda\in
\Omega$  {\sl grows with exponent\/} 
$\ge \gamma$ if there exist arbitrarily large positive
$t$ such that 
$$
\delta( g_{t} \Lambda) \le e^{-\gamma t}\,.
$$
 Also define the {\sl growth exponent\/}
$\gamma(\Lambda)$ of $\Lambda\in\Omega$ with respect to $\{g_{t}\}$
to be the supremum of all $\gamma$ for which the $g_{t}$-trajectory of
$\Lambda$  grows with exponent 
$\ge \gamma$. In view of the preceding remark, one has
$$
\gamma(\Lambda) = \limsup_{t\to\infty}
\frac{\text{dist}(g_t\Lambda, \bz^{n+1})}t\,.
$$

One sees that the orbit has a nonzero growth exponent
if and only if there is a sequence $t_j\to\infty$ such that the distance
from 

It is easy to see from the Borel-Cantelli Lemma 
and the asymptotics of 
$\nu(\ggm\sm K_\eta)$ as $\eta \to 0$ \cite[\S 7]{KM-loglaws} that $\gamma(\Lambda) = 0$
for $\nu$-a.e.~$\Lambda\in \ggm$. 
\endcomment

The passage from \da\ to growth exponents of
trajectories will be based 
on the following 
elementary lemma:

\proclaim{Lemma  2.1} Suppose we are given a set $E\subset \br^2$ which is
discrete and homogeneous with respect to positive integers, that is,
 $kE\subset E$ for any $k\in\bn$. Also take $a,b > 0$, $v > a/b$, and define $\gamma$ by 
$$
\gamma = \frac
{bv - a}{v + 1}\quad \Leftrightarrow 
\quad v = \frac{a+\gamma}{b-\gamma}\,.\tag 2.2
$$ 
Then the following are equivalent:

\roster
\item"{\rm [2.1-i]}" 
there exist $(x,z)\in E$ with
arbitrarily large
$|z|$ such that 
$|x| \le |z|^{-v}$;
\item"{\rm [2.1-ii]}" there exist arbitrarily large $t > 0$ such that for some
$(x,z)\in E\nz$ one has
$$
\max\left(e^{a t} |x|,  e^{-bt} |z|\right) \le e^{-\gamma t}\,.\tag 2.3
$$ 
\endroster
\endproclaim

\demo{Proof} Assume [2.1-i], take $(x,z)\in E$ with $|x| \le
|z|^{-v}$, and 
define $t$  by $e^{-bt} |z| = e^{-\gamma t}$, that is, $|z| = e^{(b-\gamma) t}$. (Note that it
follows from (2.2) that $\gamma < b$.) Then one has
$$
 e^{a t}|x|  \, {\le} \,
e^{a t}|z|^{-v} = e^{a t} (e^{(b-\gamma) t})^{-v}   \un{(2.2)}{=} e^{-\gamma t}
\,,
$$
that is, (2.3) holds for this choice of $x$, $z$ and $t$. Taking $|z|$  arbitrarily large
produces arbitrarily large $t$ as well.

Assume now that 
[2.1-ii] holds. Then one can find a sequence
$t_n\to\infty$ and
$(x_n,z_n)\in E\nz$ such that 
$$
e^{a t_n} |x_n| \le e^{-\gamma t_n}\quad \text{and}\quad  e^{-bt_n} |z_n|\le  e^{-\gamma t_n}
\,,\tag 2.4
$$
and write 
$$
 |x_n|  \un{(2.4)}{\le}
e^{-(a+\gamma )t_n}\un{(2.2)}{=} e^{-v(b-\gamma)t_n} \un{(2.4)}{\le}
|z_n|^{-v}\,.
$$ 
If the sequence $\{z_n\}$ is unbounded, [2.1-i] is proved. Otherwise, note that $x_n\to 0$ due to
(2.4); by the discreteness of
$E$, the sequence $\{(x_n,z_n)\}$ must stabilize, and thus one has $(0,z)\in E$ for some $z > 0$.
But then  $(0,kz)\in E$ for any $k\in\bn$ by the homogeneity, and the proof of  [2.1-i] is
finished. 
\qed\enddemo

\proclaim{Corollary  2.2} $\vy\in \Cal W_v(1,n)$ iff the growth exponent
$\gamma_{\sssize F}(u_\vy\bz^{n+1})$ of $u_\vy\bz^{n+1}$ with respect to $F$ 
as in {\rm (2.1)} is not less than $\gamma$, 
the latter being defined by 
$$
\gamma = \frac
{v-n}{n(v+1)}\quad \Leftrightarrow 
\quad v = \frac{n(1+\gamma)}{1- n \gamma }\,.\tag 2.5
$$ 
\endproclaim

\demo{Proof} This corollary is in fact a special case of 
Theorem 8.5 from \cite{KM2}. However
one can easily derive it from the previous lemma by taking $a = 1$, $b = 1/n$ and
$$
E = \big\{\big(\vy\vq + p,\|\vq\|\big)\bigm| (p,\vq)\in\bz^{n+1}\big\}\,,
$$
and noticing that the inequality
$$
\delta( g_{t} u_\vy\bz^{n+1}) \le e^{-\gamma t}\tag 2.6 
$$
amounts to the validity of (2.3) for some $(x,z)\in E\nz$.
\qed\enddemo

\proclaim{Corollary  2.3}  The following are equivalent for $\vy\in \br^n$ and $F$ 
as in {\rm (2.1)}:
\roster
\item"{\rm [2.3-i]}" $\vy$ is VWA;
\item"{\rm [2.3-ii]}" $\gamma_{\sssize F}(u_\vy\bz^{n+1}) > 0$;
\item"{\rm [2.3-iii]}" 
for some $\gamma > 0$ there exist infinitely many
$t\in\bn$ such that {\rm (2.6)} holds.
\endroster
\endproclaim

\demo{Proof} The equivalence of [2.3-i] and [2.3-ii] is straightforward 
from Corollary  2.2 and
(2.5),  while to derive [2.3-iii] one notices that the ratio of
$\delta(g_{t}\cdot)$ and $\delta(g_{t'}\cdot)$ is uniformly 
bounded from both sides when $|t-t'| < 1$.
\qed\enddemo

We return now to the setting of the \da\ on subsets of $\br^n$. More precisely, 
we consider  
a map $\vf = (f_1,\dots,f_n):U\to \br^n$, where  $U$ is an open subset of
 $\br^d$, and study \di\ properties of vectors $\vf(\vx)$ for a.e.~$\vx\in U$.  
This calls for considering the corresponding map from
$U$ into $\Omega_{n+1}$, namely
$
\vx\mapsto u_{\vf(\vx)}\bz^{n+1}$, where
$$u_{\vf(\vx)} \df \left(\matrix
1 & \vf(\vx)  \\
0 & I_n
\endmatrix \right)\,,\tag 2.7
$$
and then looking at growth
of trajectories of lattices $u_{\vf(\vx)}\bz^{n+1}$ under the action 
of $
g_t$
as in (2.1).

In the next section we will describe a method, introduced in \cite{KM1},   which
is based on keeping track on what happens to every
subgroup $\Gamma$ of $\bz^{n+1}$ under the action by $u_{\vf(\vx)}$ and then by $g_t$.
Fix a basis   
$\ve_0,\ve_1,\dots,\ve_n$ of $\br^{n+1}$, and for
$I =
\{i_1,\dots,i_j\}\subset \{0,\dots,n\}$, $i_1 < i_2 < \dots < i_j$,  let 
$$\ve_{\sssize I} \df \tsize
\ve_{i_1}\wedge\dots\wedge \ve_{i_j}\in \bigwedge^j(\br^{n+1})\,,$$
 with the convention
$\ve_\vrn = 1$. We extend the norm $\|\cdot\|$ from   $\br^{n+1}$  to the exterior algebra
$\bigwedge(\br^{n+1})$ by
$\|\sum_{I\subset \{1,\dots,j\}}w_{\sssize I}\ve_{\sssize I}\| = \max_{I\subset
\{0,\dots,n\}}|w_{\sssize I}|$. Thus it makes sense to define the norm of $\Gamma$ as above by 
$
\|\Gamma\| \df \|\vw\| $, where $\vw\text{ represents }\Gamma$.
Note that  the ratio of $
\|\Gamma\|$ and the volume of the quotient space $\br\Gamma/\Gamma$ is uniformly bounded
between two positive constants (depending on $n$ and on the choice of the norm). 
 Also note that it follows from Minkowski's
Theorem that 
$\delta(\Gamma)$ must be small whenever $
\|\Gamma\|$ is small; 
more precisely, for any $j > 0$ there exists a positive constant $c(j)$ 
such that
$$
\delta(\Gamma) \le c(j)
\|\Gamma\|^{1/j}\tag 2.8
$$
for any $\Gamma$ of rank $j$.

As a preparation for the next section, 
let us write down a formula for $g_tu_{\vf(\vx)}\vw$, where $\vw$
represents a subgroup $\Gamma$ of $\bz^{n+1}$.
\comment
We need some additional notation. 
For a {\sl lattice\/} (discrete subgroup of maximal rank) $\Lambda$  in $\br^{n+1}$
and
for $1\le j \le n$, denote by  $\Cal S_j(\Lambda)$ the set of all
subgroups of $\Lambda$ of rank $j$, and by $\Cal
S_*(\Lambda)$ the set of all
 nonzero proper  subgroups of $\Lambda$, that is, 
$\Cal
S_*(\Lambda) = \cup_{j = 1}^n \Cal
S_j(\Lambda)$. With some abuse of notation, we will also denote by $\Cal
S_j(\Lambda)$ and $\Cal
S_*(\Lambda)$ the set of $\vw\in \bigwedge(\br^{n+1})$ representing $\Gamma$ as above. 

collection of all 
subgroups of $g_{t}u_{\vf(\vx)}\bz^{n+1}$, or, equivalently, 
the
collection of elements   $g_{t}u_{\vf(\vx)}\vw$ where $\vw\in\Cal
S_*(\bz^{n+1})$. 
\endcomment
Note that the action of
$u_{\vf(\vx)}$ leaves $\ve_0$ invariant and sends $\ve_{i}$ to $\ve_{i} + f_i(\vx) \ve_0$, $i =
1,\dots,n$. Therefore 
$$
u_{\vf(\vx)}\ve_{\sssize I} =  \cases &\ve_{\sssize I} \hskip 1.9in \text{ if }0\in I\\
&\ve_{\sssize I} + \sum_{i\in I} (-1)^{l(I,i)} f_i(\vx)
\ve_{\sssize I \cup \{0\}\ssm\{i\}}\text{ otherwise}\,,\endcases
$$
where one defines 
$$l(I,i) \df \text{ the number of elements of $I$ strictly between $0$ and $i$.}
$$ 
Taking $\vw$ of the form $\sum_{I}w_{\sssize I}\ve_{\sssize I}$, 
one gets
$$
u_{\vf(\vx)}\vw =  \sum_{0\in I}\big( w_{\sssize I} + \sum_{i\notin I} (-1)^{l(I,i)}
w_{\sssize I
\cup\{i\}\ssm\{0\}}f_i(\vx) \big) \ve_{\sssize I} + \sum_{0\notin I}w_{\sssize I}\ve_{\sssize
I} \,,
$$
and, further, for $\vw\in\bigwedge^j(\br^{n+1})$,
$$
g_tu_{\vf(\vx)}\vw =  e^{\frac {n+1-j}n t} \sum_{0\in I}
\big( w_{\sssize I} + \sum_{i\notin I}
(-1)^{l(I,i)} w_{\sssize I
\cup\{i\}\ssm\{0\}}f_i(\vx) \big) \ve_{\sssize I}  + e^{-\frac jn t}
\sum_{0\notin I} w_{\sssize I}\ve_{\sssize
I}\,.\tag 2.9
$$
What is important here is that each of the coordinates of $
g_tu_{\vf(\vx)}\vw$ is expressed as a linear combination of functions
$1,f_1,\dots,f_n$.

\heading{3. Extremality criteria for  good maps}
\endheading 


Let us recall  the  definition
introduced in
\cite{KM1}. If  $C$ and $\alpha$ are  positive
numbers and $V$ a subset of $\br^d$, let us say that   a 
function
$f:V\to \br$ is {\sl \cag\ on\/}   $V$  if 
$$
\aligned
\text{for any open ball
$B\subset V$ and any $\vre > 0$, one has}\\
\big|\{x\in B\bigm| |f(x)| < \vre\cdot{\sup_{x\in
B}|f(x)|}\}\big| \le 
C\vre^\alpha |B|\,.\quad 
\endaligned\tag 3.1
$$
See \cite{KM1, BKM} for various properties and examples of \cag\ functions. One property 
will be particularly useful: it is easy to see that 
$$
\text{$f_i$, $i\in I$, are \cag\ on $V$
$\Rightarrow$ 
so is $\sup_{i\in I}|f_i|$}\,.\tag 3.2
$$

Now let $\vf = (f_1,\dots,f_n)$ be 
a  
map from an open subset
$U$ of
$\br^d$ to $\br^n$. We will say that $\vf$ is {\sl good  at  $\vx_0\in
U$\/} if there exists a neighborhood
$V \subset U$ of $\vx_0$ and positive $C,\alpha$  such that any linear combination of $
1,f_1,\dots,f_n$ is
$(C,\alpha)$-good on
$V$. 
We
will
 say that
$\vf:U\to \br^n$  is {\sl good\/} if the set of $\vx_0\in U$ such that $\vf$ is good  at 
$\vx_0$ has full measure. Note that $C,\alpha$ do not have to be uniform
in $\vx_0\in U$; however, once $V \ni \vx_0$ is chosen, every function of
the form $f = c_0 + c_1f_1 + \dots+ c_nf_n$ must satisfy (3.1) for some
uniformly chosen $C$ and $\alpha$.

Recall (see \cite{KM1, Lemma 3.2}) that the basic example of \cag\ functions is given by
polynomials: any polynomial map $\vf:\br^d\to \br^n$ is good at every
point of $\br^d$.  A more general  class of examples is given by linear
combinations of coordinate functions of nondegenerate maps:

\proclaim{Proposition 3.1 {\rm \cite{KM1, Proposition 3.4}}} Let $\vf = (f_1,\dots,f_n)$ be a
smooth map from an open subset $U$ of
$\br^d$ to $\br^n$ which is $l$-nondegenerate in $\br^n$ at  $\vx_0\in U$. Then there exists a
neighborhood
$V\subset U$ of $\vx_0$ and positive $C$  such that any 
linear combination of $ 1,f_1,\dots,f_n$
is
$(C,1/dl)$-good on
$V$. \endproclaim

In other words, $\vf$ is good at every point at which it is nondegenerate in $\br^n$. From this
one easily derives

\proclaim{Corollary  3.2} Let $\Cal L$ be an affine
subspace of $\br^n$ and let  $\vf = (f_1,\dots,f_n)$ be a
  smooth
map  from an open subset $U$ of
$\br^d$ to $\Cal L$ which is nondegenerate in $\Cal L$ at $\vx_0\in U$.
Then
$\vf$ is good at $\vx_0$.
\endproclaim

\demo{Proof} Put $\dim(\Cal L) = {s}$,  choose any
affine map $\vh$  from $\br^{s}$ onto $\Cal L$, and define 
$\vg = (g_1,\dots,g_{s})$ by $\vg = \vh^{-1}\circ
\vf$.  It follows from the nondegeneracy of $\vf$  in $\Cal L$
that $\vg$ is nondegenerate in $\br^{s}$ at $\vx_0$, hence, by Proposition
3.1,  it is good at $\vx_0$. To finish the proof it
suffices to observe that any linear combination of $
1,f_1,\dots,f_n$ is a linear combination of $
1,g_1,\dots,g_{s}$.  \qed\enddemo

\proclaim{Corollary  3.3} Let $\vf$ be a
real analytic map from a connected open subset $U$ of
$\br^d$ to $\br^n$. Then there exists an affine subspace $\Cal L$ 
of $\br^n$ such that $\vf$ is nondegenerate in $\Cal L$ at every point of
$U$; consequently, $\vf$ is good at every point of
$U$.
\endproclaim

\demo{Proof} For any $\vx\in U$, denote by $\Cal L_0(\vx)$ the linear
space spanned by all partial derivatives of $\vf$ at $\vx$, and put
$\Cal L(\vx) = \vf(\vx) + \Cal L_0(\vx)$. Then for any $\vx_0,\vx\in U$
such that  the Taylor series of $\vf$ 
centered at  
$\vx_0$ converges at $\vx$, one has $\vf(\vx)\in \Cal
L(\vx_0)$ and $\Cal L(\vx) \subset \Cal L(\vx_0)$. Since $U$ is
connected, for any 
$\vx'\in U$ one can find a finite sequence $\vx_1,\dots,\vx_k = \vx'$ such
that the Taylor series of $\vf$ centered at $\vx_{i-1}$ converges at
$\vx_{i}$ for all $i = 1,\dots,k$. Therefore $\Cal
L(\vx') = \Cal L(\vx_0)$, and, by reversing the roles of 
$\vx_0$ and $\vx'$, one sees that $\Cal
L \df \Cal L(\vx)$ is independent
of
$\vx\in U$. It remains to notice that from the construction it follows
that
$\vf$ is nondegenerate in
$\Cal L$ at every 
$\vx\in U$, and apply Corollary 3.2. 
 \qed\enddemo

Note that it follows from the proof that $\Cal L$ can be defined as the 
intersection of all the affine subspaces of $\br^n$ containing $\vf(U)$, or, equivalently,
as $$\vf(\vx_0) + \text{Span}\,\{\vf(\vx) - \vf(\vx_0)\mid \vx\in U\}$$ for any $\vx_0\in U$.

\example{Example 3.4} It is instructive for better understanding of the
class of good maps to remark that the assumption of the analyticity of
$\vf$ cannot be dropped. Indeed, let us sketch a construction of a
$C^\infty$ function from
$[0,1]$ to
$\br_{\sssize +}$ which is not good on a subset of $[0,1]$ of positive measure.
First  for every $J = (a,b)\subset [0,1]$ define 
$$
\psi_{\sssize J}(x) \df \varphi(x-a) \varphi(b-x)\,,\quad\text{where }\varphi(x) \df
\cases 0\,,\hskip .35in x\le 0\,, \\ e^{-1/x^2},\ x \ge 0\,.\endcases
$$
One can easily  verify that for every 
neighborhood $V$ of either
$a$ or $b$ it is impossible to find
$C,\alpha > 0$  such that $\psi_{\sssize J}$ is \cag\ on $V$. 
Then consider a Cantor set $K\subset[0,1]$ of positive measure, and for
$k\in\bn$ let
$\Cal J_k$ be the collection of disjoint subsegments of $[0,1]$
thrown away at the $k$th stage of the construction of $K$. (For example,
one can divide every interval left at the $k$th stage onto $3^{k+1}$ equal
pieces and then throw away the middle interval.)
After that define 
$$
\psi(x) \df \sum_{k = 1}^\infty c_k\sum_{J\in\Cal J_k}\psi_{\sssize J}(x)\,,
$$
where $c_k$ decays fast enough as $k\to\infty$ to guarantee that $\psi$
is $C^\infty$ (in the aforementioned example, one can take $c_k = 3^{-3^k}$). Since every
neighborhood of every point of
$K$  contains an endpoint of $J\in\Cal J_k$ for some $k\in\bn$, it follows
that
$\psi$ is not good at any $x\in K$. 
\endexample

We now state  an  estimate from \cite{KM1}, which will be used to derive 
a criterion for
the extremality of  $\vf$ once the latter is
chosen within the class of good maps. It will be convenient to use the following notation: if $B =
B(\vx,r)$ is a ball in $\br^d$ and $c > 0$, we will denote by $cB$ the
ball
$B(\vx,cr)$.

 \proclaim{Theorem 3.5 {\rm (cf.~\cite{KM1, Theorem 5.2})}} For any  $d,k\in\bn$ there exists
a positive constant $C'$ (explicitly 
estimated
in \cite{KM1}) such
that the following holds. Given 
$C,\alpha > 0$,
$0 <
\rho  \le 1/k$,   a ball $B\subset \br^d$ and a continuous map $h:3^k B
\to
\GL_k(\br)$, 
let us assume that 
\roster
\item"{\rm [3.5-i]}" for any
$\Gamma \in \Cal S_k$, the function $x\mapsto \|h(\vx)\Gamma\|$ is 
\cag\ on $3^kB$, and
\item"{\rm [3.5-ii]}" for any
$\Gamma \in \Cal S_k$,  $\sup_{\vx\in B}\|h(\vx)\Gamma\|\ge \rho $.
\endroster
Then 
 for any  positive $ \vre \le \rho$ one has
$$
\left|\{\vx\in B\mid \delta\big(h(\vx)\bz^k\big) < \vre\}\right| \le C C'
\left(\frac\vre \rho \right)^\alpha  |B|\,.
$$
\endproclaim

Informally speaking, the conclusion of the above theorem  says that 
the  `orbit' 
$\{ h(\vx)\bz^k \mid \vx \in B \}\subset \Omega_k$ `does not
diverge', that is, its very significant 
proportion (computed in terms of Lebesgue measure on $B$)  stays inside
compact sets \linebreak
$\{\lio_k\mid \delta(\Lambda) \ge \vre\}$. We remark that such
nondivergence
results have a long history, dating back to the 
work of Margulis \cite{Ma} in the 1970s, and many applications in the
theory of
dynamics on homogeneous spaces, see e.g.~\cite{KSS, Chapter 3}
 for a historical account.

\medskip

The next 
lemma  sharpens  \cite{KM1, Theorem 5.4}, giving a condition sufficient
for the extremality of  a good map $\vf$.


 \proclaim{Lemma 3.6}  Let $B$ be a ball in $\br^d$, and let $\vf =
(f_1,\dots,f_n)$ be a continuous map from $3^{n+1}B$  to $\br^n$. 
Suppose that:

\roster
\item"{\rm [3.6-i]}" $\exists\,C,\alpha > 0$ such that any linear
combination of
$ 1,f_1,\dots,f_n$ is
$(C,\alpha)$-good on
$3^{n+1}
B$;
\item"{\rm [3.6-ii]}" for any $\beta
> 0$  there exists $T = T(\beta) > 0$ such that for any $t\ge T$ and
any $\Gamma \in \Cal S_{n+1}$ one has
$$
\sup_{\vx \in B} \|g_t u_{\vf(\vx)}\Gamma\|\ge e^{-\beta t}\,,
$$
where $u_{\vf(\vx)}$ is as
in {\rm (2.7)} and $F = \{g_t\}$ is as in {\rm (2.1)}.
\endroster 
Then $\vf(\vx)$ is not VWA for a.e.\ $\vx\in B$.
\endproclaim

\demo{Proof} We apply Theorem 3.5 with
$k = n+1$ and 
$h(\vx) = g_t u_{\vf(\vx)}$.  Our goal is
 to
show that for any $\gamma > 0$, the set 
$
\big\{\vx\in B\mid \gamma_{\sssize F}\big(u_{\vf(\vx)}\bz^{n+1}\big) > \gamma\big\}
$
has measure zero. 
As was observed  in the preceding section,  see (2.9), for
every
$\vw
\in
\Cal S_{n+1}$  all the coordinates of $h(\vx)\vw$ are  linear combinations of
$1,f_1,\dots,f_n$, which, by (3.2), implies [3.5-i].  Now choose any
$\beta < \gamma$, take $t \ge \max\big((T(\beta), \log(n+1)/\beta\big)$
and put $\rho = e^{-\beta t}$. This guarantees $\rho  \le 1/k$ and 
verifies condition [3.5-ii]. Thus, for
$t$ as above, the measure of each of the sets
$$\{\vx\in B\mid \delta(g_{t}u_{\vf(\vx)}\bz^{n+1})
< e^{-\gamma t}\}\tag 3.3
$$
is not greater than $CC'e^{-\alpha(\gamma-\beta) t}$. One then 
applies 
the Borel-Cantelli Lemma to conclude that almost every
$\vx\in B$ belongs to at most  finitely many of the sets (3.3) with $t\in\bn$,
which completes the proof in view of   Corollary 2.3.
\qed\enddemo

The next lemma shows that assumption [3.6-ii] is in fact necessary
for the extremality of 
$\vf$; furthermore, the consequences of  [3.6-ii] being not true
are much stronger than positive measure of the set $\{\vx\in
B\mid\vf(\vx)\text { is  VWA}\}$.

 \proclaim{Lemma 3.7}  Let $B$ be a ball in $\br^d$, and let $\vf$ be a 
 map from $
B$  to $\br^n$ such that {\rm [3.6-ii]} does not hold. Then $\vf(\vx)$ is
VWA for all $\vx\in B$. 
\endproclaim

\demo{Proof} The assumption says that 
there exists   $\beta > 0$ such that 
one has
$$\sup_{\vx \in B}\|g_t u_{\vf(\vx)}\Gamma\| < e^{-\beta t}$$ for
arbitrarily large $t$ (and $\Gamma\in \Cal S_{n+1}$ dependent on $t$). In
view of (2.8), for any $\vx \in B$ this implies
$$\delta(g_t u_{\vf(\vx)}\bz^{n+1}) \le c(j)\|g_t u_{\vf(\vx)}\Gamma\|^{1/j}
< c(j)e^{-\beta t/j}\,,$$ 
where $j$ is the rank of $\Gamma$.  Hence
$\gamma_{\sssize F}\big(u_{\vf(\vx)}\bz^{n+1}\big) \ge \beta/n$, and an application
of 
Corollary 2.3 finishes the proof. 
\qed\enddemo

We now combine the two lemmas above to obtain the desired extremality
criterion.

 \proclaim{Theorem 3.8} Let $U$  be an  open subset $\br^d$, and let 
$\vf$ be a map from $U$  to $\br^n$ which is continuous and good. 
Then the following are equivalent:

\roster
\item"{\rm [3.8-i]}" the set $\{\vx\in
U\mid\vf(\vx)\text { is not VWA}\}$ is dense in $U$;
\item"{\rm [3.8-ii]}" $\vf$ is extremal (that is, the above set
has full measure);
\item"{\rm [3.8-iii]}" for a.e.~$\vx_0\in U$ and any $\,r > 0$ there exists
a ball
$B\subset U$ centered at
$\vx_0$ of radius less than $r$ satisfying {\rm [3.6-ii]};
\item"{\rm [3.8-iv]}" any ball $B\subset U$ satisfies {\rm [3.6-ii]}.
\endroster 
\endproclaim

\demo{Proof} Obviously [3.8-ii]$\Rightarrow$[3.8-i] and
[3.8-iv]$\Rightarrow$[3.8-iii]. The implication \linebreak
[3.8-i]$\Rightarrow$[3.8-iv] is immediate from Lemma 3.7. 
Assuming [3.8-iii] and using the fact that $\vf$ is good, for
a.e.~$\vx_0\in U$ one finds a ball $B$ centered at
$\vx_0$ such that the dilated ball $3^{n+1}B$ is contained in $U$,
and both [3.6-i] and [3.6-ii] hold. Thus Lemma 3.6 applies, and  [3.8-ii]
follows.
\qed\enddemo

We remark that for the equivalence of [3.8-i] and [3.8-ii] it
is essential that $\vf$ is
chosen within the class of good maps.
Indeed, one can consider the map $\vf(x) =  \big(x,\psi(x)\big)$, with
$\psi$ from Example 3.4.  Clearly $\vf(x)$ is VWA for any $x$ from $K$,
which is assumed to have positive measure.
On the other hand, the restriction of $\vf$ to any $J\in\Cal J_k$, $k\in
\bn$, is nondegenerate in $\br^2$, hence the set  $\{\vx\in
(0,1)\mid\vf(\vx)\text { is not VWA}\}$ has full measure in $(0,1)\ssm
K$, and the latter is dense in
$(0,1)$.

\medskip

Our next task is to 
rephrase [3.6-ii].
For any $B\subset U$ let us
denote by 
$\Cal F_B$ the 
$\br$-linear span of the restrictions of $1,f_1,\dots,f_n$ to $B$. Then, for   any
ball $B\subset U$,  let us assume that  functions
$g_1,\dots,g_{s}: B\to \br$ are chosen so that
$1,g_1,\dots,g_{s}$ form a basis of $\Cal F_{B}$ (here the dimension ${s}+1$ of $\Cal F_B$ may
depend on $B$).    Further,  the choice of functions
$1,g_1,\dots,g_{s}$ defines 
a matrix $P\in M_{{s}+1,n+1}$ formed by
coefficients in the expansion of
$1,f_{1},\dots,f_n$ as linear combinations of $1,g_1,\dots,g_{s}$. 
In other words, with the notation $\tilde \vf \df (1,f_1,\dots,f_n)$ and 
$\tilde \vg \df (1,g_1,\dots,g_{s})$,
one has 
$$ 
\tilde \vf(\vx) = \tilde
\vg(\vx)  P\quad\forall\,\vx\in B\,.\tag 3.4
$$
 This way, the
restriction of
$\vf$  to $B$ is described by two pieces of data: the $({s}+1)$-tuple 
$\tilde
\vg$ and the  matrix $P$. We now proceed to show that, assuming the map $\vf$ is good
(which is an assumption involving $\tilde
\vg$),
a criterion for its extremality can be written in terms of  $P$.

Indeed,  any $f\in\Cal F_B$ can be written as $f(\vx) = \tilde \vg(\vx) \vv$ for some
$\vv\in\br^{n+1}$, and because  of the linear independence of components
of
$\tilde \vg$ over $\br$, the `supremum-on-$B$' norm of $f$, that is,
$f\mapsto \sup_{\vx\in B}|f(\vx)|$, is equivalent to
$\|\vv\|$, the constant in the equivalence depending  on $B$ and the choice of  $\tilde
\vg$. Now recall that for any $\vw = \sum_{I}w_{\sssize I}\ve_{\sssize I}$, the $I$th component
of  $u_{\vf(\vx)}\vw$ is equal to $
 w_{\sssize I}$ if ${0\notin I}$, and to 
$$
w_{\sssize I} +
\sum_{i\notin I} (-1)^{l(I,i)} w_{\sssize I
\cup\{i\}\ssm\{0\}}f_i(\vx) 
$$ 
if $I$ contains $0$. It will be convenient to simplify the
latter expression by introducing the following notation:
given  $I\subset
\{0,\dots,n\}$ containing $0$ with $|I| = j$, and an element $\vw = \sum_{I}w_{\sssize
I}\ve_{\sssize I}$ of $\bigwedge^j(\br^{n+1})$, 
let us define a vector $\vc_{\sssize I,\vw}\in\br^{n+1}$
by 
$$
\vc_{\sssize I,\vw} \df \sum_{i\notin (I\nz)} (-1)^{l(I,i)} w_{\sssize I
\cup\{i\}\ssm\{0\}}\ve_i \ \  = \ \ w_{\sssize
I}\ve_0 + \sum_{i\notin I} (-1)^{l(I,i)} w_{\sssize I
\cup\{i\}\ssm\{0\}}\ve_i\,.\tag 3.5
$$ 
Then the nonconstant components of $u_{\vf(\vx)}\vw$ can be written as $\,\tilde\vf(\vx)
\vc_{\sssize I,\vw} = \tilde
\vg(\vx)  P\vc_{\sssize I,\vw}$; therefore the `supremum-on-$B$' norm of each of the functions
can be replaced by $\| P\vc_{\sssize I,\vw}\|$.  Modifying all the norms and using (2.9), one 
replaces 
$
\,\sup_{\vx\in B}\|g_tu_{\vf(\vx)}\vw\|$ by  $\max\big(e^{\frac {n+1-j}n t} \max_{0\in I} \|P
\vc_{\sssize I,\vw}\|, e^{-\frac jn t} \max_{0\notin I} |w_{\sssize
I}| \big)$.
We summarize the above discussion as follows:

 \proclaim{Proposition 3.9} Let $B$ be a ball in $\br^d$, $\vf$  a 
 map from $B$  to $\br^n$, $\{1,g_1,\dots,g_{s}\}$ a basis  of
$\Cal F_B$,  and $P$ a matrix 
 satisfying  {\rm (3.4)}.   Then {\rm [3.6-ii]} is  equivalent to
$$
\aligned
\forall\,\beta
> 0\ \ \exists\, T  > 0\ \text{ such that }\ \forall \,t\ge T\,, \ \forall\,j = 1,\dots,n
 \text{ and
} \forall\,\vw \in \Cal S&_{n+1,j}\\ \text{one has}\quad
\max\big(e^{\frac {n+1-j}n t} \max_{0\in I} \|P
\vc_{\sssize I,\vw}\|, e^{-\frac jn t} \max_{0\notin I} |w_{\sssize
I}| \big)\ge e^{-\beta t}\,.
\endaligned\tag 3.6
$$
\endproclaim

\comment
\item"{\rm [3.5-ii]}" $\forall\,k =
1,\dots,n$ and $\ \forall\,v > \frac {n+1-k}k$
$\  \exists\, N  > 0$ such that  for any $\ \vw \in \Cal S_{n+1,k}$ \newline with 
$\ \max_{0\notin I}
|w_{\sssize I}| > N$, one has $$\max_{0\in I} \|P
\vc_{\sssize I,\vw}\| > (\max_{0\notin I}
|w_{\sssize I}|)^{-v}\,.$$

\demo{Proof} The preceding argument shows the equivalence [3.6-ii]
$\Leftrightarrow$ [3.5-i], while 
the fact that [3.5-i] and
[3.5-ii] are equivalent is a special case of Lemma 2.1,  where for each $k =
1,\dots,n$ one considers
$$E = \big\{\big(\max_{0\in I} \|P
\vc_{\sssize I,\vw}\|,\max_{0\notin I}
|w_{\sssize I}|\big) \bigm|\vw \in \Cal S_{n+1,k}\big\}\,.\quad\qed$$
\enddemo
\endcomment

In other words, we have shown that the extremality of a continuous  good map $\vf:U\to\br^n$ is 
 equivalent to the validity of 
certain \di\  conditions
involving  matrices $P$ `coordinatizing' $\vf|_B$. These conditions will be made more precise in
the next section, and Theorem 1.2 will be obtained as a corollary.

\heading{4. Extremality criteria for  affine subspaces}
\endheading

Note that 
in general, as was mentioned above, it may not be possible to choose
the same matrix $P$ uniformly for all balls $B$ in $U$. Let us now consider an important 
special case when
this is  possible:

 \proclaim{Theorem 4.1} Let $U\subset \br^d$ be an open subset,  and let
$\vg = (g_1,\dots,g_{s}): U\to \br^s$ be a continuous good map  
such that  
$$ 
\aligned
\forall\,\vx\in U\ &\text{(equivalently, $\forall\,\vx$ from a dense subset of $U$)}\\
\text{the germs of }&1,g_1,\dots,g_{s}\text{ 
at $\vx$ are linearly independent over $\br$.}
\endaligned\tag 4.1
$$
Also fix $P\in M_{{s}+1,n+1}$, and let $\vf$ be given by {\rm (3.4)}. 
Then  {\rm (3.6)} is equivalent to each of the following conditions:

\roster
\item"{\rm [4.1-i]}" the set $\{\vx\in
U\mid\vf(\vx)\text { is not VWA}\}$ is non-empty;
\item"{\rm [4.1-ii]}" $\vf$ is extremal (that is, the above set
has full measure).
\endroster 
\endproclaim

\demo{Proof} It is clear from (3.4) that $\vf$ is also continuous and good.
[4.1-ii] follows from
(3.6) in view of Proposition 3.9 and the implication [3.8-iii]$\Rightarrow$[3.8-ii]. 
On the other hand, if (3.6) is violated, Lemmas 3.7 and Proposition 3.9
imply that  every point of $U$
has a neighborhood $B$ such that $\vf(\vx)$ is
VWA for all $\vx\in B$, contradicting [4.1-i].
\qed\enddemo

Informally speaking, the assumption of Theorem 4.1 says that $\vf$ is not `assembled
from several pieces unrelated to each other'. Without that assumption, one can
easily construct examples of good continuous (and even $C^\infty$) 
maps 
$\vf$ 
satisfying [4.1-i] but not [4.1-ii].

\medskip

Now suppose that  $\vf$ is real analytic and $U$ is connected, or, more generally, 
that $\vf$ is  nondegenerate 
in some affine subspace 
$\Cal L$ of
$\br^n$. Then one can easily 
find 
$s \le n$ and a
good $s$-tuple  
$\vg$ satisfying
 (4.1). Specifically, one takes $s = \dim(\Cal L)$ and, as in the proof of Corollary 3.2, 
defines $\vg$ to be equal to  $\vh^{-1}\circ
\vf$ where $\vh$ is any
affine map   from $\br^{s}$ onto $\Cal L$. Furthermore, one easily recovers the 
corresponding matrix $P$ by  writing $\vh$ in the form $$\tilde
\vh(\vx) =  \tilde \vx P\,,\tag 4.2
$$
where as usually we have $\tilde
\vh \df (1,h_1,\dots,h_n)$. 
It follows that for fixed $\Cal L$,
any $\vf$  which is nondegenerate in  $\Cal L$ will satisfy the 
assumptions of Theorem 4.1 with  some uniformly chosen  $P$. In particular, Theorem
4.1 applies to the map $\vf = \vh$ given by (4.2), that is, to the
subspace $\Cal L$ itself. Thus we have proved:

 \proclaim{Theorem 4.2} Let $\Cal L$ be an $s$-dimensional affine subspace  of $\br^n$
 parametrized 
 as
in {\rm (4.2)} with 
$P\in M_{s+1,n+1}$. Then each of the following conditions below is 
equivalent to {\rm (3.6)}:

\roster
\item"{\rm [4.2-i]}" $\Cal L$ contains at least one not \vwa\ point;
\item"{\rm [4.2-ii]}" $\Cal L$ is extremal;
\item"{\rm [4.2-iii]}" any smooth submanifold of $\Cal L$ which is nondegenerate in $\Cal L$ is
extremal.
\endroster 
\endproclaim

One then recovers Theorem 1.2 as the implications
[4.2-ii]$\Rightarrow$[4.2-iii] and \linebreak [4.2-i]$\Rightarrow$[4.2-ii]
 above. Note also that Theorem 1.1 is obtained by taking $s = n$ and $P = I_{n+1}$.

\medskip

Now recall that, as described in the introduction, one can by permuting variables  
without loss of generality parametrize $\Cal L$ by (1.4)
 for some $A\in M_{s+1,{n-s}}$; that is, take $P$ of the form $P = \pmatrix I_{s+1} &
A\endpmatrix$. In order to restate Theorem 4.2 in terms of
$A$, let us denote by $\vc_{\sssize I,\vw}^{\sssize +}$ (resp.~$\vc_{\sssize I,\vw}^{\sssize -}$) the column vector
consisting of  the first
$s+1$ (resp.~the last $n-s$) coordinates of $\vc_{\sssize I,\vw}$. In other words, we  have 
$\vc_{\sssize I,\vw} = \pmatrix \vc_{\sssize I,\vw}^{\sssize +} \\ \vc_{\sssize I,\vw}^{\sssize -}\endpmatrix$ 
where
$$
\vc_{\sssize I,\vw}^{\sssize +} = \sum_{i\in \{0\}\cup(\{1,\dots,s\}\ssm I)} (-1)^{l(I,i)} w_{\sssize I
\cup\{i\}\ssm\{0\}}\ve_i  =  w_{\sssize
I}\ve_0 + \sum_{i\in \{1,\dots,s\}\ssm I} (-1)^{l(I,i)} w_{\sssize I
\cup\{i\}\ssm\{0\}}\ve_i
$$ 
and
$$
\vc_{\sssize I,\vw}^{\sssize -} = \sum_{i\in \{s+1,\dots,n\}\ssm I} (-1)^{l(I,i)} w_{\sssize I
\cup\{i\}\ssm\{0\}}\ve_i\,.
$$ 
Then one has:

 \proclaim{Theorem 4.3}  Let $\Cal L$ be an $s$-dimensional affine subspace  of $\br^n$
 parametrized 
 as
in {\rm (1.4)} with  
$A\in M_{s+1,n-s}$. Then the following are
equivalent:
\roster
\item"{\rm [4.3-i]}" $\Cal L$ is extremal ($\Leftrightarrow$ {\rm [4.2-i]} 
and {\rm [4.2-iii]}
hold);
\item"{\rm [4.3-ii]}" for any
$\beta > 0$  there exists $T  > 0$ such that  for any $t\ge T$,  $j =
1,\dots,n$ and  $\vw \in \Cal S_{n+1,j}$ one has
$$
\max\big(e^{\frac {n+1-j}n t} \max_{0\in I} \|
\vc_{\sssize I,\vw}^{\sssize +} + A\vc_{\sssize I,\vw}^{\sssize -}\|, e^{-\frac jn t} \max_{0\notin I} |w_{\sssize
I}| \big)\ge e^{-\beta t}\,;
$$
\item"{\rm [4.3-iii]}" $\forall\,j =
1,\dots,n\text{ and }\ \forall\,v > \tfrac {n+1-j}j
\quad  \exists\, N  > 0$  such that  for any $\ \vw \in \Cal S_{n+1,j}${  with }
\ $\max_{0\notin I}
|w_{\sssize I}| > N\,,$ one has 
$$\max_{0\in I} \|\vc_{\sssize I,\vw}^{\sssize +} + A\vc_{\sssize
I,\vw}^{\sssize -}\| > (\max_{0\notin I} |w_{\sssize I}|)^{-v}\,.$$
\endroster 
\endproclaim
 
\demo{Proof} The preceding argument shows the equivalence [4.3-i]
$\Leftrightarrow$ [4.3-ii], while 
the fact that [4.3-ii] and
[4.3-iii] are equivalent is a special case of Lemma 2.1,  where for each $j =
1,\dots,n$ one considers
$$E = \big\{\big(\max_{0\in I} \|\vc_{\sssize I,\vw}^{\sssize +} + A\vc_{\sssize I,\vw}^{\sssize -}\|,\max_{0\notin I}
|w_{\sssize I}|\big) \bigm|\vw \in \Cal S_{n+1,j}\big\}\,.\quad\qed$$
\enddemo
 
It is instructive to write down a special case of the above inequality corresponding to $j = 1$. 
That is, let us take $\vv \in  \bz^{n+1}\nz = \Cal S_{n+1,1}$ in place of $\vw$; one sees
that the only one-element subset $I$ of $\{0,\dots,n\}$ for which $\vc_{\sssize I,\vv}$
is defined is $I = \{0\}$, and it easily follows from (3.5) that $\vc_{\sssize \{0\},\vv} =
\vv$. Writing 
$\vv  = \pmatrix 
\vp
\\
\vq
\endpmatrix$, where  $\vp\in\bz^{s+1}$ and $\vq\in\bz^{n-s}$, one gets $\vc_{\sssize \{0\},\vv}^{\sssize +}
= \vp$ and $\vc_{\sssize \{0\},\vv}^{\sssize -}
= \vq$. Then denoting by $\vp'$ the vector with components $p_1,\dots,p_s$, one writes the
$j = 1$ case of [4.3-iii]  as follows:

\roster
\item"{\rm [4.3-iii]$_{j = 1}$}" {\it for any $\,v > n$
there exist at most finitely many $\vq\in\bz^{n-s}$ such that for some  $\vp\in\bz^{s+1}$  one
has }
$$ \|\vp + A\vq\|
\le \left\|\pmatrix 
\vp'
\\
\vq
\endpmatrix\right\|^{-v}\,.\tag 4.3
$$
\endroster 

Now observe that one can safely replace the latter inequality by 
$$\|\vp + A\vq\|
\le \|
\vq
\|^{-v}\,,\tag 4.4
$$ perhaps slightly changing $v$. Indeed, (4.4) clearly follows from (4.3). On the other hand,
(4.4) implies that $\|\vp\| \le C\|\vq\|$ for some
$C$ dependent only on $A$; thus, for a slightly smaller $v$ and large enough $\|\vq\|$, (4.3)
would follow. We arrive to the conclusion that  
[4.3-iii]$_{j = 1}$ is 
equivalent to 
(1.7).

However, let us point out that one does not need the full strength of Theorem 4.3 to see that
(1.7) is one of the conditions necessary for the extremality of $\Cal L$  
as in {\rm (1.4)}. Indeed, as shown above, the assumption  $A\in\Cal W_n^{\sssize +}({s}+1,n-{s})$ 
amounts to the existence of  $v > n$
such that for  infinitely many
$\vq\in
\bz^{n-{s}}$ one can find $\vp\in\bz^{{s}+1}$ satisfying (4.3). Then one can take  any
$\vx\in \br^{s}$  and write
$$
\big|\,p_0 + (\vx,
\tilde\vx A) \left(\matrix \vp' \\ \vq
\endmatrix\right)\big|  = |p_0 + \vx\vp' + \tilde\vx A\vq| = |\tilde\vx( A\vq +
\vp)|\le ({s}+1)\left\|\tilde\vx\right\|
\| A\vq + \vp\|\,.
$$ 
Slightly decreasing $v$ if needed, one gets infinitely many
solutions of 
$$
\left|\,p_0 + (\vx,
\tilde\vx  A) \left(\matrix \vp' \\ \vq
\endmatrix\right)\right|  \le \left\|\pmatrix 
\vp'
\\
\vq
\endpmatrix\right\|^{-v}\,,
$$
that is,  $(\vx,
\tilde\vx \tilde A)$ is proved to be VWA for all $\vx$. 

\medskip

Let us now ask the following 
question: {\it could it be the case that 
the remaining
$n-1$ conditions of Theorem 4.3 are redundant, that is,  follow from\/} [4.3-iii]$_{j = 1}$? 
The affirmative answer to this question would provide
a very easy to state extremality criterion, i.e.~the validity of (1.7), for affine subspaces
and their submanifolds.

The answer to this question is (in general) not known to the author. However, the next result
shows that the case $j = n$ of [4.3-iii] is indeed redundant.

 \proclaim{Lemma 4.5} For any $\,{s} = 1,\dots,n-1$, any $\,A\in M_{s+1,n-s}$ and any $\,\vw \in
\Cal S_{n+1,n}$ one has
$$\max_{0\in I} \|\vc_{\sssize I,\vw}^{\sssize +} + A\vc_{\sssize I,\vw}^{\sssize -}\|
\ge 1\,.\tag 4.5
$$
\endproclaim

\demo{Proof} Denote by $J_i$ the set $\{0,\dots,n\}\ssm\{i\}$. Then any $\vw \in \Cal
S_{n+1,n}$ can be written in the form
$\vw = \sum_{i = 0}^n w_i \ve_{\sssize J_i}$, and from (3.5) it follows that for $i =
1,\dots,n$ one has
$$\vc_{\sssize J_i,\vw} = w_i\ve_0 + (-1)^{i-1}w_0\ve_i\,.\tag 4.6
$$
Therefore for any $i = 1,\dots,
 s$ one has  $\vc_{\sssize J_i,\vw}^{\sssize -} = 0$, and  hence 
$$
\vc_{\sssize J_i,\vw}^{\sssize +} + A\vc_{\sssize J_i,\vw}^{\sssize -} = \vc_{\sssize J_i,\vw}^{\sssize +}  = w_i\ve_0 +
(-1)^{i-1}w_0\ve_i\,.
$$
Consequently, (4.5) is satisfied whenever $w_0 \ne 0$. On the other
hand, $w_0 = 0$, in view of (4.6), implies that $\vc_{\sssize J_i,\vw}^{\sssize -} = 0$ for any $i$. 
Taking $i > s$ for which $w_i \ne 0$, one gets 
$
\|\vc_{\sssize J_i,\vw}^{\sssize +} + A\vc_{\sssize J_i,\vw}^{\sssize -}\| =  |w_i|.
$ \qed\enddemo

This, in particular, gives an affirmative answer to the above question in the case $n = 2$:  a
line in
$\br^2$ given by 
$y = a_0 + a_1 x$ is extremal if and only if  $\left(\matrix a_0 \\ a_1
\endmatrix\right)\notin \Cal W_2^{\sssize +}(2,1) $.

It turns out that an   argument similar to the proof of Lemma 4.5  produces an analogous
extremality criterion for  $(n-1)$-dimensional affine subspaces of $\br^n$ for arbitrary $n$:

 \proclaim{Lemma 4.6} Let  $\,A\in M_{n,1}$ be given by a column vector $\va\in\br^{n}$ (this
corresponds to
$s = n-1$). Then
{\rm (4.5)} holds for  any  
$\,\vw
\in
\Cal S_{n+1,j}$ with
$j > 1$. 
\endproclaim

\demo{Proof} Our choice of $s$ implies that for any $\,\vw
\in
\Cal S_{n+1,j}$ and any $I\ni 0$ of size $j$, the vector $\vc_{\sssize I,\vw}^{\sssize -}$ consists of
 a single
component, namely
$$
c_{\sssize I,\vw}^{\sssize -} = \cases (-1)^{l(I,n)} w_{\sssize I
\cup\{n\}\ssm\{0\}} \text{ if } n\notin I\\ 0 \hskip 1.1in\text{ otherwise.}
\endcases\tag 4.7
$$ Therefore for $0,n\in I$ one can
write
$$
\vc_{\sssize I,\vw}^{\sssize +} + \va c_{\sssize I,\vw}^{\sssize -} = \vc_{\sssize I,\vw}^{\sssize +}  = w_{\sssize
I}\ve_0 + \sum_{i\in\{1,\dots,n-1\}\ssm I} (-1)^{l(I,i)} w_{\sssize I
\cup\{i\}\ssm\{0\}}\ve_i\,.
$$
Consequently, (4.5) is satisfied whenever $w_{\sssize
I} \ne 0$ for some $I$ containing  $n$. (Here we use the fact that $j > 1$: indeed, such an $I$
must also contain some $i = 0,\dots,n-1$, and thus the $i$th coordinate of $\vc_{\sssize
I\cup\{0\}\ssm\{i\},\vw}^{\sssize +}$ is equal to $w_{\sssize
I}$.) On the other hand, the assumption 
$w_{\sssize
I} = 0$ for all $I\ni n$, in view of (4.7), implies that
$c_{\sssize I,\vw}^{\sssize -} = 0$ for any
$I$; thus (4.5) is satisfied again. Hence one can take an arbitrary  $I$ for which $w_{\sssize
I} \ne 0$  and observe that for any $i\in I$, the absolute value of the 
$i$th coordinate of
$\vc_{\sssize I\cup\{0\}\ssm\{i\},\vw}^{\sssize +}$ is equal to $|w_{\sssize
I}|$.
\qed\enddemo

Combining the above lemma with Theorem 4.3 and the equivalence of [4.3-iii]$_{j
= 1}$ and  (1.7), one  easily obtains Theorem 1.3 under assumption (1.5). In particular,
in view of (1.8), one sees that the set of $\va\in\br^{n}$ for which the  $(n-1)$-dimensional
affine subspace  of $\br^n$ parametrized  by
$$
\vx = (x_1,\dots,x_{n-1})\mapsto(\vx,
\tilde\vx \va) = (x_1,\dots,x_{n-1},a_0 + a_1x_1+\dots+a_{n-1}x_{n-1})\tag 4.8
$$
is not extremal, has \hd\ $1$.

\medskip 

As was mentioned above, we are unable to show the equivalence of
(1.7) and [4.3-iii] in the general case. However, let us now turn to
assumption (1.6), under which that equivalence can be demonstrated by a
direct proof (that is, without a reference to lattices).

Until the end of this section, let us assume $s = 1$, take $A$  of the form
$\left(\matrix 0 \\ \vb
\endmatrix\right)$ for a row vector $\vb= (b_1,\dots,b_{n-1})\in\br^{n-1}$, and let $\Cal L$ be
parametrized by (1.4); that is, $\Cal L$ is a line passing through the origin given by 
$$
x\mapsto(x,b_1x,\dots,b_{n-1}x)\,.\tag 4.9
$$
 It is
clear that 
$A\in \Cal W_n^{\sssize +}(2,n-1)$ if and only if $\vb\in \Cal W_n^{\sssize +}(1,n-1)$. 
Therefore in order to prove Theorem 1.3 assuming (1.6), it suffices to prove the
following

 \proclaim{Proposition 4.7} A line given by {\rm (4.9)} is extremal whenever 
$$
\vb\notin \Cal
W_n^{\sssize +}(1,n-1)\,.\tag 4.10
$$
\endproclaim

\demo{Proof} We follow an argument from the paper \cite{BBDD}, where a stronger (than
the extremality) property was proved for $\Cal L$ as in (4.9) under the stronger (than (4.10)) 
assumption   ${\vb\notin \Cal
W_n^{\sssize -}(1,n-1)}$.

The goal is to prove that for any  $v>n$, the set
$$
\left\{x\in \br\left| 
\aligned
|p + q_0x + q_1b_1x + \dots + q_{n-1}b_{n-1}x| \le \|\vq\|^{-v} \\
\text{ for
infinitely many }\vq = (q_0,q_1, \dots, q_{n-1})^{\sssize T} \in\bz^n, \ &p\in\bz
\endaligned
\right.\right\}\tag 4.11
$$
has measure zero. Clearly without loss of generality one can restrict $x$ to lie in the unit
interval. Also,  our usual  notation $\tilde \vb = (1,\vb)$ will be helpful, since the left
hand side of the inequality in (4.11) will be then written as $|p + (\tilde \vb \vq) x|$.

Let us now
state a lemma from which the desired result will easily follow.
For $\vb\in\br^{n-1}$ and $Q,v > 0$, define the set
$\Cal A(\vb,v,Q)$ to be the set of $x\in [0,1)$ for which the inequality
$$
|p + (\tilde \vb \vq) x| < Q^{-v}\tag 4.12
$$
holds for some $p\in\bz, \vq\in\bz^n$ with $Q\le \|\vq\| < 2Q$. 

 \proclaim{Lemma 4.8} For any  $\vb$ satisfying {\rm (4.10)} and any $v > n$ there
exists a positive constant $C = C(\vb,v)$ such that for any $Q > 1$ one has
$$
|\Cal A(\vb,v,Q)| < CQ^{\frac{n-v}2}\,.
$$
\endproclaim

It is easy to see that the intersection of the set (4.11) with $[0,1)$ is contained in 
$$
\left\{x \mid  x\in \Cal A(\vb,v,2^k)\text{ for infinitely many }k\in\bn\right\}\,.\tag 4.13
$$
Assuming Lemma 4.8, one has  $|\Cal A(\vb,v,2^k)| < C2^{-\frac{v-n}2k}$ $\forall\,k$, and
the fact that the set (4.13) has measure zero is then immediate from the Borel-Cantelli
Lemma. \qed\enddemo

It remains to write down the

\demo{Proof of Lemma 4.8}  Define $\Cal A^0(\vb,v,2^k)$ to be the set of $x\in [0,1)$ for which (4.12)
holds for some $\vq\in\bz^n$ with $Q\le \|\vq\| < 2Q$ and with $p = 0$. It is contained in a
union of intervals of the form $\big[0,\frac {Q^{-v}}{|\tilde \vb \vq|}\big]$.  Due to (4.10),
there exists $c = c(v) > 0$ such that the denominator of the above fraction is not less than 
$$c\cdot \max(q_1,\dots,q_{n-1})^{-\frac{n+v}2} \ge c\cdot \|\vq\|^{-\frac{n+v}2} <
c\cdot (2Q)^{-\frac{n+v}2}\,.
$$
Therefore one has
$$
|\Cal A^0(\vb,v,Q)| < c^{-1} Q^{-v} (2Q)^{\frac{n+v}2} = c^{-1} 2^{\frac{n+v}2}
Q^{\frac{n-v}2}\,.
$$

Now let us estimate the measure of $\Cal A(\vb,v,Q) \ssm \Cal A^0(\vb,v,Q)$. Note that, assuming $p\ne 0$,
inequality (4.12) can be solvable in $x\in [0,1)$ only if $|\tilde \vb \vq| > 1- Q^{-v} > 1 -
1/Q \ge 1/2$. For fixed $p$ and $\vq$, (4.12) defines an interval of length at most  $2\frac
{Q^{-v}}{|\tilde \vb \vq|}$, and, for fixed $\vq$, the number of
different centers of those intervals, that is, points $\frac {p}{|\tilde \vb \vq|}$, is at most
$1 + |\tilde \vb \vq|$. Therefore one can write
$$
\split
|\Cal A(\vb,v,Q) \ssm \Cal A^0(\vb,v,Q)|\ &\le \sum\Sb\|\vq\|< 2Q,\\|\tilde \vb \vq| > 1/2\endSb  \frac
{2Q^{-v}}{|\tilde \vb \vq|} (1 + |\tilde \vb \vq|)   = \ 2Q^{-v} \sum\Sb\|\vq\|< 2Q,\\|\tilde
\vb \vq| > 1/2\endSb \left(1 + \frac {1 }{|\tilde \vb \vq|}\right) \\  &\le\  2Q^{-v} (4Q)^n +
2Q^{-v}\sum\Sb\|\vq\|< 2Q,\\|\tilde \vb \vq| > 1/2\endSb \frac
{1 }{|\tilde \vb \vq|}\,.
\endsplit
$$
To estimate the sum in the right hand side of the above formula, note that for fixed
$q_1,\dots,q_{n-1}$ and variable $q_0$, the values of $\tilde \vb \vq$ form an arithmetic
progression. Thus, fixing $q_1,\dots,q_{n-1}$, one gets
$$\split
\sum\Sb\|q_0\|< 2Q,\\|\tilde \vb \vq| > 1/2\endSb \frac
{1 }{|\tilde \vb \vq|} &\le 2\left(\frac
{1 }{1/2} + \frac
{1 }{1/2 + 1} + \dots +\frac{1}{1/2 + 2Q-1}\right) \\ &= 4\left(1 + \frac13 + \dots +
\frac1{4Q-1}\right) < 4\big(1 + \log (4Q-1)\big) 
\,.
\endsplit$$

Summing the above estimate over all $q_1,\dots,q_{n-1}$, one obtains
$$
|\Cal A(\vb,v,Q) \ssm \Cal A^0(\vb,v,Q)|\ \le 2^{2n+1}Q^{n-v} + 8 Q^{-v} (4Q)^{n-1}\big(1 + \log
(4Q-1)\big)\,,
$$
which is not greater than the right hand side of the desired inequality for an appropriate value
of
$C$ .
\qed\enddemo

\heading{5. Multiplicative approximation}
\endheading 

The dynamical approach to \di\ problems described above has an advantage of being quite
general to allow various modifications of the set-up. In particular, most of the ideas described
in this paper work for the so-called \ma. Let us briefly list all the relevant definitions.
For $\vx\in\br^n$ we let 
$$
\Pi_{\sssize +}(\vx) 
= \prod_{i = 1}^n |x_i|_{\sssize +}\,,\quad\text{where}\quad|x|_{\sssize +}=\max(|x|, 1)\,.
$$
For $v > 0$ let us denote by $\Cal {WM}_v(1,n)$ the set of row vectors
$\vy\in\br^n$ for which  there are infinitely
many
$\vq\in
\bz^n$ such that
$$
 |\vy\vq + p|   \le \Pi_{\sssize +}(\vq)^{-v/n} \quad \text{for some
}p\in\bz\,.  \tag 5.1
$$
Since the right hand side of (5.1) is not less than that of (1.1), one clearly has
$\Cal {WM}_v(1,n) \supset \Cal
W_{v}(1,n)$; in particular,  $\Cal {WM}_{n}(1,n) = \br^n$ by Dirichlet's
Theorem. Also it can be shown using  the Borel-Cantelli Lemma that the Lebesgue measure of $\Cal
{WM}_v(1,n)$ is zero whenever
$ v > n$. Therefore, with the definition of   {\sl very well
multiplicatively approximable\/}  (VWMA) vectors as those $ \vy\in\br^n$ which are in  $\Cal
{WM}_v(1,n)$ for some
$ v > n$, one has that almost all $ \vy\in\br^n$ are not VWMA.

Let us now say, 
following the 
terminology of \cite{Sp4},  that a  submanifold 
 $\Cal M$ of $\br^n$ (resp.~a smooth map $\vf$ from an open subset $U$ of $\br^d$ to $\br^n$) is 
 {\sl
strongly extremal\/} if almost all $\vy\in \Cal M$ (resp.~$\vf(\vx)$  for a.e.~$\vx\in
U$) are not
VWMA. It is clear that strong extremality implies extremality, and to prove a manifold to be
strongly extremal is usually a harder task than just to prove extremality. For example, the 
strong extremality of the curve (1.2) (that is, the multiplicative analogue of Mahler's
problem) was conjectured by A.~Baker in 1975 \cite{B}, and the only 
proof that exists as of now is based on the dynamical approach of  \cite{KM1}.
In fact, the main result of the latter paper (\cite{KM1, Theorem A}, of which Theorem 1.1 is a
special case)  is  the strong extremality of  manifolds nondegenerate in
$\br^n$ (in the analytic case this was conjectured in \cite{Sp4}). 

With the help of the approach developed in \cite{KM1}, let us now try to investigate \ma\
properties of  generic points of  proper affine subspaces and their submanifolds by first
describing
the set of VWMA vectors in a dynamical language. It turns out that the actions that are
relevant for this case are multi-parameter. Namely, one replaces (2.1) by 
$$
\aligned
&g_\vt = \text{\rm diag}(e^{t},e^{-t_1},\dots,e^{-t_n})\,,\\
\text{ where }\ \  &\vt = (t_1,\dots,t_n),
t_i \ge 0, \text{ and  }t = \sum_{i = 1}^nt_i\,.
\endaligned\tag 5.2
$$
The latter notation is used throughout the section, so that whenever $t$ and $\vt$ 
appear in the same context, $t$  stands for $\sum_{i = 1}^nt_i$.

For the rest of this section, we  mostly sketch our argument, as it is
very similar to what is done in \S\S 2--4, and  only highlight important
modifications. The following is a multi-parameter version of Lemma 2.1:

\proclaim{Lemma  5.1} Suppose we are given a set $E$ of pairs $(x,\vz)\in\br^{n+1}$, which is
discrete and homogeneous with respect to positive integers. 
 Then
the following are equivalent:

\roster
\item"{\rm [5.1-i]}" 
for any $\,v > n$ there exist  $(x,\vz)\in E$ 
with $\vz$ arbitrarily far from $0$
 such that 
$$|x| \le \Pi_{\sssize +}(\vz)^{-v/n}\,;\tag 5.3
$$
\item"{\rm [5.2-ii]}" for any $\,\gamma > 0$ there exists an unbounded set of $\vt\in
\br^{n}_{\sssize +}$ 
for which one has 
$$
e^{ t} |x| \le e^{-\gamma t}\quad \text{and}\quad  e^{-t_i} |z_i|\le  e^{-\gamma t},\quad i =
1,\dots,n
\,.\tag 5.4
$$
for some
$(x,\vz)\in E\nz$.
\endroster
\endproclaim

\demo{Proof} We follow the lines of the proof of Lemma 2.1. Assuming [5.1-i], take $v >
n$ and $(x,\vz)\in E$ 
satisfying (5.3), and define $t$  by 
$$e^{(1-n\gamma) t} = \Pi_{\sssize +}(\vz)\,,\tag 5.5a$$
where $\gamma < 1/n$ is as in (2.5). Then for every $i$ define $t_i$  by
$$e^{t_i} = e^{\gamma t}|z_i|_{\sssize +}\,.\tag 5.5b$$ 
Note that, since $|z_i| \le |z_i|_{\sssize +}$, this implies $e^{-t_i}
|z_i|\le e^{-\gamma t}$; and  note also that multiplying all the
equalities (5.5b) and comparing the result with (5.5a) one can verify that $t = \sum_{i =
1}^nt_i$.  Then  one has
$$
 e^{t} |x|  \, {\le} \,
e^{ t}\Pi_{\sssize +}(\vz)^{-v/n} = e^{ t} (e^{(1-n\gamma) t})^{-v/n}   \un{(2.5)}{=} e^{-\gamma
t}
\,,
$$
that is, (5.4) is satisfied for this choice of $x$, $\vz$ and $\vt$; taking $\vz$ with 
arbitrarily large $\Pi_{\sssize +}(\vz)$ produces arbitrarily large values of $t$.

Assume now that 
[5.1-ii] holds, and take $\gamma < 1/n$. Then one can find an unbounded  sequence of vectors
$\vt$  and a sequence of points 
$(x,\vz)\in E\nz$ satisfying (5.4). Since for any $\vt$ one has $t_i \ge \gamma t $ for at least
one $i$,  passing to a subsequence and reshuffling the coordinates of $\vt$  and $\vz$ if
necessary, one can assume that  for some $k = 1,\dots,n$ and all entries $\vt$ of that sequence,
one has
$$
t_i \ge \gamma t \quad \text{for}\quad i \le k, \qquad \text{and}\quad t_i < \gamma t \quad
\text{for}\quad i > k\,.\tag 5.6
$$
It follows from (5.6) and (5.4) that $|z_i|_{\sssize +} \le e^{t_i - \gamma t}$ for $i \le k$,
and  $|z_i| < 1$ for
$i > k$, hence
$$
\Pi_{\sssize +}(\vz) \le
\prod_{i = 1}^k |z_i|_{\sssize +}  \le  e^{t_1 + \dots + t_k - k\gamma t} \le e^{(1 - k\gamma)
t}\,.\tag 5.7$$
Now it is time to find an appropriate $v$. However, because of  an extra
parameter
$k$, we have to modify (2.5), namely  define $v > n$ by
$$
\gamma = \frac
{v-n}{ kv + n}\quad \Leftrightarrow 
\quad v = \frac{n(1+\gamma)}{1- k \gamma }\,.\tag 5.8
$$
 Then the right hand side of (5.7) is equal to $e^{\frac nv(1+\gamma)t}$, and (5.4) implies
$$
| x | \le e^{-(1+\gamma)  t} = (e^{(1 - k\gamma)
t})^{-v/n}\un{(5.7)}{\le}\Pi_{\sssize +}(\vz)^{-v/n}\,.
$$
which is exactly what was needed. After that, as in the proof of Lemma 2.1, one notices that
a uniform bound on $\|\vz\|$,  
by the discreteness of
$E$, implies that $(0,\vz_0)\in E$ for some $\vz_0$, and  integral multiples of  $(0,\vz_0)$
give infinitely many $(x,\vz)\in E$  satisfying (5.3). To finish the proof, it remains  to
observe that (5.8) forces
$v$ to tend to $n$ uniformly in $k$ as $\gamma \to 0$.
\qed\enddemo

\proclaim{Corollary  5.2}  For $\vy\in
\br^n$ and $g_\vt$ as in {\rm (5.2)}, the following are equivalent:

\roster
\item"{\rm [5.2-i]}" $\vy$ is VWMA;
\item"{\rm [5.2-ii]}" 
for some $\gamma > 0$ there exists an unbounded set of $\vt\in
\br^{n}_{\sssize +}$ 
such that 
$$\delta(g_{\vt}u_\vy\bz^{n+1}) \le e^{-\gamma t}\,;\tag 5.9
$$
\item"{\rm [5.2-iii]}" 
for some $\gamma > 0$ there exist infinitely many
$\vt\in\bz_{\sssize +}^n$ such that {\rm (5.9)} holds.
\endroster
\endproclaim

Note that Corollary 2.2 in \cite{KM1} provides the
implication [5.2-i]$\Rightarrow$[5.2-iii].

\demo{Proof} Taking
$
E = u_\vy\bz^{n+1}$, 
one sees that (5.9) amounts to the validity of (5.4) for some $(x,\vz)\in E\nz$. The
rest of the argument mimics the proof of Corollary 2.3.
\qed\enddemo

From the above corollary and Theorem 3.5 it is not hard to derive multiplicative analogues of 
extremality criteria of \S 3 and \S4. The crucial condition to consider is an analogue of  
[3.6-ii]: if $B\subset \br^d$ is a ball and $\vf$ a map from $B$ to $\br^n$, it is important 
to check whether or not 
$$
\aligned
\forall\,\beta
> 0\ \ \exists\, T= T(\beta)  > 0\ \text{such that  for any }\vt\in
\br^{n}_{\sssize +} \text{ with }t&\ge  T\\
\text{ and any }\Gamma \in \Cal S_{n+1}, \text{ one has}\quad
\sup_{\vx \in B} \|g_t u_{\vf(\vx)}\Gamma\|\ge e^{-\beta t}\,.
\endaligned\tag 5.10
$$
The following can be proved by a straightforward repetition of the argument of \S 3:

\roster
\item"$\bullet$" (cf.\ Lemma 3.7) if {\rm (5.10)} does not hold, then $\vf(\vx)$ is
VWMA for all $\vx\in B$;
\item"$\bullet$" (cf.\ Lemma 3.6) if $\vf$ is continuous, defined on $3^{n+1}B$ and satisfies
{\rm [3.6-i]} and
{\rm (5.10)}, then $\vf(\vx)$ is not VWMA for a.e.\ $\vx\in B$.
\endroster

Therefore one has

 \proclaim{Theorem 5.3} For $U$   and 
$\vf$  as in Theorem 3.8, the following are equivalent:
\roster
\item"{\rm [5.3-i]}" the set $\{\vx\in
U\mid\vf(\vx)\text { is not VWMA}\}$ is dense in $U$;
\item"{\rm [5.3-ii]}" $\vf$ is strongly extremal (that is, the above set
has full measure);
\item"{\rm [5.3-iii]}" for a.e.~$\vx_0\in U$ and any $\,r > 0$ there exists
a ball
$B\subset U$ centered at
$\vx_0$ of radius less than $r$ satisfying {\rm (5.10)};
\item"{\rm [5.3-iv]}" any ball $B\subset U$ satisfies {\rm (5.10)}.
\endroster 
\endproclaim

In order to  express (5.10) in \di\ language,
we need  some more notation: for $\vt\in\br_{\sssize +}^n$ and $I\subset
\{0,\dots,n\}$, let
$$
t_{\sssize I} \df \sum_{i\in I\nz}t_i\,.
$$
Then, taking $\vf$, $B$, $\{1,g_1,\dots,g_{s}\}$  and $P$ as in Proposition 3.9, one can observe that
 (5.10) can be written in the form
$$
\aligned
\forall\,\beta
> 0\ \ \exists\, T  > 0\ \text{ such that }\ \forall \,\vt\in
\br^{n}_{\sssize +} \text{ with }t\ge  T  \text{ and
} \forall\,\vw \in \Cal S&_{n+1}\\
\text{ one has}\quad
\max\big(e^{t - t_{\sssize I}} \max_{0\in I} \|P
\vc_{\sssize I,\vw}\|, e^{- t_{\sssize I}} \max_{0\notin I} |w_{\sssize
I}| \big)\ge e^{-\beta t}\,.
\endaligned\tag 5.11
$$ 

Thus Theorem 5.3 implies

 \proclaim{Theorem 5.4} Let $U$, $\vf$, 
$\vg$ and  $P$ be as in Theorem 4.1. 
Then  {\rm (5.11)} is equivalent to each of the following conditions:

\roster
\item"{\rm [5.4-i]}" the set $\{\vx\in
U\mid\vf(\vx)\text { is not VWMA}\}$ is non-empty;
\item"{\rm [5.4-ii]}" $\vf$ is strongly extremal (that is, the above set
has full measure).
\endroster 
\endproclaim

Taking $P$ of the form $P = \pmatrix I_{s+1} &
A\endpmatrix$, one deduces

 \proclaim{Theorem 5.5} The
following are equivalent for an $s$-dimensional affine subspace $\Cal L$ of $\br^n$ parametrized 
as in {\rm (1.4)} with 
$A\in M_{s+1,n-s}\,$:

\roster
\item"{\rm [5.5-i]}" $\Cal L$ contains at least one not VWMA point;
\item"{\rm [5.5-ii]}" $\Cal L$ is  strongly extremal;
\item"{\rm [5.5-iii]}" any smooth submanifold of $\Cal L$ which is nondegenerate in $\Cal L$ is
strongly  extremal;
\item"{\rm [5.5-iv]}" for any
$\beta > 0$  there exists $T  > 0$ such that  for any $\vt\in\br_{\sssize +}^n$ with $t\ge T$,  
 and any $\vw \in \Cal S_{n+1}$ one has
$$
\max\big(\max_{0\in I} e^{t - t_{\sssize I}} \|
\vc_{\sssize I,\vw}^{\sssize +} + A\vc_{\sssize I,\vw}^{\sssize -}\|, \max_{0\notin I} e^{- t_{\sssize I}}
|w_{\sssize I}| \big)\ge e^{-\beta t}\,.\tag 5.12
$$
\endroster 
\endproclaim
  
This, in particular, proves Theorem 1.4, as well as gives a 
criterion for
the strong extremality of  $\Cal L$ written in terms of \di\ properties of 
the parametrizing matrix
$A$. 

\medskip

In general, it seems to be a hard problem to simplify condition [5.5-iv]. However, in view of
computations made in
\S 4, this can be easily done in the case $s = n-1$, that is, when $\Cal L$ is a codimension one
subspace. Namely, due to Lemma 4.6, one knows that  (5.12) always holds when $s =
n-1$ and $\vw \in \Cal S_{n+1,j}$ with $j > 1$; thus it suffices to handle the
case $j = 1$. As in \S 4, one then uses  $\vv  = \pmatrix 
\vp
\\
q
\endpmatrix\in  \bz^{n+1}\nz 
$ in place
of
$\vw$ (here $\vp\in\bz^{n}$ and $q\in\bz$) and  notices that   $\vc_{\sssize
I,\vv}$ is defined only for  $I = \{0\}$, with 
$\vc_{\sssize
\{0\},\vv} =
\vv$, $\vc_{\sssize \{0\},\vv}^{\sssize +}
= \vp$ and $\vc_{\sssize \{0\},\vv}^{\sssize -}
= q$. Replacing $A$ by a column vector $\va = (a_0,a_1,\dots,a_{n-1})^{\sssize T}\in\br^{n-1}$
and letting $\vp= (p_0,p_1,\dots,p_{n-1})^{\sssize T}$ and $\vp'=
(p_1,\dots,p_{n-1})^{\sssize T}$, one gets 

 \proclaim{Corollary 5.6} The
following are equivalent for an $(n-1)$-dimensional affine subspace $\Cal L$ of $\br^n$ parametrized 
by {\rm (4.8)}:

\roster
\item"{\rm [5.6-i]}" $\Cal L$ is  strongly extremal; 
\item"{\rm [5.6-ii]}" for any
$\beta > 0$  there exists $T  > 0$ such that  for any $\vt\in\br_{\sssize +}^n$ with $t\ge T$ and any 
$\pmatrix 
\vp
\\
q
\endpmatrix\in  \bz^{n+1}\nz$ one has
$$
\max\big( e^{t} \|
\vp + \va q \|,  e^{- t_1}|p_1|,\dots,e^{- t_{n-1}}|p_{n-1}|, e^{- t_{n}}|q|
 \big)\ge e^{-\beta t}\,.
$$
\item"{\rm [5.6-iii]}" for any $v > n$
there exists $K  > 0$ such that for any $\vp\in  \bz^{n}$ and $q\in  \bz$ with
$\,\max(\|\vp'\|,|q|)
 > K$, one has 
$$ \|
\vp + \va q \|
> \Pi_{\sssize +}(
\vp',
q
)^{-v/n}\,.$$\endroster 
\endproclaim

\demo{Proof} The only part that requires a comment is the equivalence
[5.6-ii]$\Leftrightarrow$[5.6-iii], which is a special case of Lemma 5.1 with 
$$E = \{ (\|
\vp + \va q \|,\vp',q)\mid \vp\in\bz^{n},q\in\bz\}\,. \quad \qed
$$
\enddemo

 \proclaim{Corollary 5.7} Let $\Cal L$ be parametrized 
by {\rm (4.8)}, and let 
$$k = \#\{1\le i \le n - 1 \mid a_i \ne 0\}\,.$$ Then 
$\Cal L$ is strongly extremal iff $\,\va\notin
\Cal W_{k+1}^{\sssize +}(n,{1})$.
\endproclaim
  
\demo{Proof} By the previous corollary, 
the fact that $\Cal L$ is not strongly  extremal  is equivalent  to saying that
for some $v > n$ there exist $(\vp,q)\in\bz^{n}$ with arbitrarily large $\,\max(\|\vp'\|,|q|)
$ such that 
$$ 
\|
\vp + \va q \|
\le \Pi_{\sssize +}(
\vp',
q
)^{-v/n}\,.\tag 5.13
$$ 
Equivalently, there exists a sequence of solutions of (5.13) with $p_i$ arbitrarily
close to $a_iq$ for all
$i= 1,\dots,n-1$, which in particular happens if and only if $p_i$  is equal to zero for any
$i$ with $a_i = 0$.
Hence the ratio of $\Pi_{\sssize +}(
\vp',
q
)$ and $|q|^{k+1}$ is bounded \linebreak  from both sides, and, slightly changing $v$, one gets
infinitely many solutions of \linebreak $\|
\vp + \va q \|
\le |q|^{-\frac vn(k+1)}
$. 
\qed
\enddemo

In particular, one can see that the set of vectors $\va\in\br^{n}$ for which the  map (4.8)
is not strongly extremal is slightly bigger than the one corresponding to  non-extremality,
agrees with the latter outside of all the coordinate planes, and still has
\hd\
$1$.

\heading{6. Further generalizations and open questions}
\endheading

\subhead{6.1}\endsubhead As was mentioned before, it would be very interesting to find out
whether Theorem 1.3 can be extended to the cases when  the rank of $A$ is greater than one. 
 If it cannot,
it would be nice to find a reasonable description of the set of non-extremal
subspaces, e.g.~such that would allow to compute its \hd. Similar questions
are open in the case of multiplicative approximation; in particular, it is not clear, except
for the case of hyperplanes, how much smaller than the set of extremal subspaces is the set of
strongly extremal ones.

\subhead{6.2}\endsubhead Given that one of the main results of this paper is that the
extremality of an affine subspace is inherited by its nondegenerate submanifolds, one can ask
whether any `passage of information' takes place when the subspace is not extremal. Let us
introduce the following definition: for \amr, define the {\sl \de\/} $\omega(A)$ of $A$ by
$$
\omega(A) \df \sup\{v\mid A\in \Cal W_v(m,n)\}\,.
$$
Clearly $n/m \le \omega(A) \le \infty$ for all $A$, and $A$ is VWA iff $\omega(A) > n/m$. Now,
for a map $\vf:U\to \br^n$ define $\omega(\vf)$ to be the essential infimum of
$\omega\big(\vf(\cdot)\big)$, i.e.
$$
\omega(\vf)\df \sup\left\{v\Bigm| \big|\{\vx\in U\mid \vf(\vx)\in \Cal W_v(1,n)\}\big| >
0\right\}\,.
$$
Naturally, if $\Cal M$ is a smooth manifold, we let the \de\ $\omega(\Cal M)$ of $\Cal M$ to be the \de\ of its
parametrizing map.

Clearly a manifold (map) is extremal iff its \de\ is the smallest possible, i.e.~is equal to $n
= \omega(\br^n)$. Now the following two questions arise:

\roster
\item"$\bullet$" is it true that the \de\ of an affine subspace $\Cal L$ of $\br^n$ is inherited by
manifolds nondegenerate in $\Cal L$? 
\item"$\bullet$" how to efficiently describe the class of affine subspaces with a given \de?
\endroster

The answer to the first question is `yes', which can be proved by a refinement of the
`dynamical' approach developed in \cite{KM1} and the present paper. Furthermore, similarly to
Theorems  3.8 and 4.2--4.3, for any $v\ge n$ one can write down necessary and sufficient
conditions for a good (resp.~affine) map $\vf$ to have $\omega(\vf)\le v$. 

The answer to the second question is as
obscure as the problem of extending Theorem 1.3 beyond the  cases (1.5) and
(1.6). Indeed, the conclusion of the latter  theorem amounts to saying that, for
$\Cal L$ parametrized as in (1.4),
$\omega(A)
\le n$ implies
$\omega(\Cal L) = n$. In general, one can easily prove that
$\omega(\Cal L) \ge \omega(A)$; thus one is left to ask whether $\omega(\Cal L)$ is always
equal to $\max\big(n,\omega(A)\big)$ (this can be verified in the `rank-one' cases (1.5) and
(1.6)). All this is going to be the topic of a forthcoming paper \cite{K2}. 

\subhead{6.3}\endsubhead Even more generally, one can replace the right hand side of (1.1) by an
arbitrary function of
$\|\vq\|$.
Let us specialize to the case of row
vectors and introduce the following definition: for a non-increasing function
$\psi:\bn\to (0,\infty)$, define  $\Cal W_\psi(1,n)$ to be the set of 
$\vy\in\br^n$ for which  there are infinitely
many
$\vq\in
\bz^n$ such that
$$
 \|\vy\vq + p\|   \le \psi(\|\vq\|) \quad \text{for some
}p\in\bz\,. 
$$
It is a theorem of A.\,V.~Groshev (\cite{G}, see also \cite{S1}) that almost
no (resp.~almost all)
$\vy\in\br^n$ belong to
$\Cal W_\psi(1,n)$ if the series
$$
\sum_{k = 1}^\infty {k^{n-1}\psi(k)}\tag 6.1
$$
converges (resp.~diverges). (The
case of convergence easily follows from the Borel-Cantelli Lemma.) It has been recently
proved, see \cite{BKM} for the convergence part and \cite{BBKM} for the divergence part, that the
same dichotomy takes place for any nondegenerate submanifold of
$\br^n$;  in other words, for a smooth map $\vf:U\to \br^n$  which is nondegenerate in
$\br^n$, one has 
$$
\aligned
\text{$\vf(\vx)\in\Cal W_\psi(1,n)$ for almost no
(resp.~almost all) }&\vx\in U\\
\text{ if the series
(6.1)
converges (resp.~diverges).}
\endaligned\tag 6.2
$$ It is natural to expect  a similar dichotomy  for any
smooth submanifold  $\Cal M$ of
$\br^n$, with the convergence/divergence of (6.1) replaced by another `dividing line' condition,
possibly involving the \de\ of $\Cal M$. The following problems remain open:

\roster
\item"$\bullet$" is it true that the aforementioned `dividing line' condition of an affine
subspace
$\Cal L$ of
$\br^n$ is inherited by manifolds nondegenerate in $\Cal L$? 
\item"$\bullet$" given an affine subspace, find its `dividing line' condition; or, vice versa, 
describe the class of  subspaces with a given `dividing line'.
\endroster

The only result along these lines   known to the author is the paper  \cite{BBDD}, where it
is shown that the convergence/divergence of (6.1) serves as the `dividing line' condition for 
one-dimensional subspaces  of
$\br^n$ of the form (1.6) with ${\vb\notin \Cal
W_n^{\sssize -}(1,n-1)}$.

Finally, let us note   that the paper \cite{BKM} also contains a more general (in
particular, multiplicative) version of the convergence case of (6.2), and it is of considerable
interest to see if the argument from that paper
can be applied to the set-up of proper affine subspaces and their nondegenerate submanifolds.

\heading{ Acknowledgements}\endheading
 The author is  grateful to Professor Gregory Margulis for conjecturing
some of the results of this paper and for many stimulating and enlightening discussions, and to
Victor Beresnevich,   Barak Weiss and the reviewer   for
useful remarks.


\bigskip
\Refs

\widestnumber\key{BBDD}


\ref\key {B}\by A. Baker  \book Transcendental number theory 
\publ Cambridge Univ. Press \publaddr Cambridge
\yr 1975 \endref



\ref\key BBDD \by V. Beresnevich, V. Bernik, H. Dickinson, and M.\,M. Dodson
\paper On linear manifolds for which the Khinchin approximation theorem holds \jour Vestsi Nats.
  Acad. Navuk Belarusi. Ser. Fiz.-Mat. Navuk  \yr 2000 \pages 14--17 
\lang Belorussian \endref


\ref\key BBKM \by V. Beresnevich, V. Bernik, D. Kleinbock, and G.\,A.
Margulis 
 \paper Metric Diophantine approximation: the Khintchine--Groshev theorem
for non-degenerate manifolds \jour Moscow Math. J. \yr
2002 \vol 2 \issue 2 \pages 203--225\endref


\ref\key BD \by V. Bernik and
M.\,M. Dodson \book Metric \da\
on manifolds \publ Cambridge Univ. Press \publaddr Cambridge
\yr 1999 \endref



\ref\key BKM \by V. Bernik, D. Kleinbock, and G.\,A. Margulis \paper
Khintchine-type theorems  on
manifolds:  the convergence case for standard  and multiplicative
  versions \jour Internat. Math. Res. Notices \yr 2001   
\pages 453--486 \issue 9
\endref


\ref\key C \by J. W. S. Cassels \book An introduction to \di\ approximation \bookinfo Cambridge Tracts in Math. \vol 45
\publ Cambridge Univ. Press \publaddr Cambridge
\yr 1957 \endref
 



\ref\key {D}\by S.\,G. Dani \paper Divergent trajectories of flows on
\hs s and Diophantine approximation\jour
J. Reine Angew. Math.\vol 359\pages 55--89\yr 1985\endref

\ref\key Do  \by M.\,M. Dodson \paper 
Hausdorff dimension, lower order and Khintchine's theorem in metric Diophantine approximation
\jour  
J. Reine Angew. Math. \vol 432 \yr 1992 \pages 69--76\endref
 


\ref\key G \by A.\,V. Groshev  \paper Une th\'eor\`eme sur les
syst\`emes des formes lin\'eaires \jour Dokl. Akad. Nauk SSSR \vol 9
\yr 1938 \pages 151--152 \endref  



\ref\key {K1} \by D. Kleinbock \paper Some applications of
homogeneous dynamics to number theory \inbook in: Smooth Ergodic Theory
and Its Applications (Seattle, WA, 1999) \pages 639--660 \bookinfo Proc.
Symp. Pure Math. \vol 68  \publ Amer. Math. Soc. \publaddr Providence, RI
\yr 2001 \endref

\ref\key {K2} \bysame  \paper  An extension of quantitative nondivergence 
and applications to 
Diophantine exponents \paperinfo in preparation \endref


\ref\key {KLW} \by D. Kleinbock, E. Lindenstrauss, and B. Weiss \paper On fractal
measures and Diophantine approximation \paperinfo Preprint \yr 2003\endref

\ref\key KM1 \by D. Kleinbock and G.\,A. Margulis \paper Flows  on
homogeneous spaces and \da\ on manifolds\jour Ann. Math. \vol 148 \yr
1998 \pages 339--360 
 \endref

\ref\key KM2 \bysame \paper Logarithm laws for flows  on
homogeneous spaces \jour Invent. Math.\vol 138 \pages 451--494 \yr 1999 
\endref

\ref\key {KSS} \by D. Kleinbock, N. Shah, and A. Starkov \paper
Dynamics of subgroup actions on homogeneous spaces of Lie groups and 
applications to number theory\inbook in: Handbook on
Dynamical Systems, Volume 1A  \publ
Elsevier Science \publaddr North Holland \pages  813--930 \yr 2002  \endref

\ref\key {M}\by K. Mahler 
\paper \" Uber das Mass der Menge aller $S$-Zahlen \jour Math. Ann. \vol 106 \pages 131--139 \yr 1932\endref

\ref\key {Ma}\by G.\,A. Margulis
\paper On the action of unipotent group in the space of lattices 
\inbook Proceedings of the Summer School on group representations, (Budapest 1971)\pages
365--370\publ Acad\'emiai Kiado
\publaddr Budapest \yr 1975\endref

\ref \key R \by M.\,S. Raghunathan \book Discrete subgroups of Lie
groups 
\publ Springer-Verlag \publaddr Berlin and New York \yr 1972 \endref%
 
\ref\key {S1}\by W. Schmidt \paper  A metrical theorem in \da \jour Canadian J. Math. \vol 12
\pages 619--631 \yr 1960 \endref



\ref\key {S2}\bysame  \paper Metrische S\"atze \"uber simultane 
Approximation abh\"anginger Gr\"ossen \linebreak \jour Monatsch. Math. \vol 68 \pages
154--166\yr 1964\endref


\ref\key {S3}\bysame \paper Diophantine approximation and certain
sequences of 
lattices \jour Acta Arith. \vol 18 \yr 1971 \pages 195--178\endref

 

\ref\key {Si}\by  C.\,L. Siegel \book Zur Reduktionstheorie
quadratischer Formen \publ  The Mathematical Society of Japan \publaddr
Tokyo \yr 1959 \endref



\ref\key {Sp1}\by V. Sprind\v zuk \paper More on Mahler's conjecture  \jour  Doklady  Akad.  Nauk  SSSR \vol 155 \yr 1964 \pages 54--56  \lang Russian  \transl\nofrills English transl. in  \jour
Soviet Math. Dokl \vol 5 \pages
361--363\yr 1964\endref



\ref\key {Sp2}\bysame \book Mahler's problem in metric number theory \bookinfo Translations of Mathematical
Monographs, vol. 25 \publ Amer. Math. Soc.\publaddr Providence, RI \yr 1969 \endref

\ref\key {Sp3}\bysame \book Metric theory of Diophantine
approximations \publ
John Wiley \& Sons \publaddr New York-Toronto-London \yr 1979\endref

\ref\key {Sp4}\bysame  \paper Achievements and problems in
Diophantine approximation theory \jour Russian Math. Surveys  \vol 35 \yr 1980 \pages 1--80 \endref



\ref\key {W} \by B. Weiss \paper  Dynamics on parameter
spaces: submanifold and fractal subset questions  \inbook In:
Rigidity in Dynamics and Geometry, M. Burger and A. Iozzi
(eds.) \publ Springer \yr 2002\endref  
\endRefs
\enddocument